\documentclass[a4paper,12pt]{amsart}

\usepackage{hyperref}

\usepackage[top=3cm,left=2.5cm,right=2.5cm,bottom=3cm]{geometry}
\usepackage{relsize}
\usepackage{lineno}

\usepackage{amsmath,amsthm,pb-diagram,amssymb,comment}
\usepackage{amsfonts,graphicx,color}
\usepackage{enumitem, fancyhdr, dsfont, multicol}
\usepackage{thmtools,cleveref}
\usepackage[normalem]{ulem}
\usepackage{stmaryrd}
\usepackage{textcomp}
\usepackage{mathtools}
\usepackage{fontawesome}
\usepackage{textcomp}
\usepackage[all]{xy}

\usepackage{tikz,float}
\usepackage{mleftright}
\usepackage{mathrsfs}
\usepackage{tikz-cd}
\usepackage{dsfont}
\usepackage{stmaryrd}

\usepackage{caption}
\usepackage{subcaption}

\hypersetup{
    colorlinks=true,
    linkcolor=magenta,
    filecolor=magenta,
    urlcolor=magenta,
    citecolor=magenta,
}
\definecolor{caribbeangreen}{rgb}{0.0, 0.8, 0.6}

\setlist{topsep=0ex,itemsep=1ex}

    \DeclareMathOperator{\dom}{dom}

    \DeclareMathOperator{\ran}{ran}

    \newcommand{\menos}{\smallsetminus}

   \DeclareMathOperator{\supsum}{\overline{\rm{S}}}
    \DeclareMathOperator{\infsum}{\underline{\rm{S}}}

\newcommand{\Pbf}{\mathbf{P}}

\newcommand{\ZFC}{\mathrm{ZFC}}

\newcommand{\bfm}{\mathbf{m}}

\newcommand{\bfP}{\mathbf{P}}

\newcommand{\bbN}{\mathbb{N}}

\newcommand{\bbQ}{\mathbb{Q}}
\newcommand{\bbR}{\mathbb{R}}

\newcommand{\cB}{\mathscr{B}}
\newcommand{\cR}{\mathscr{R}}
\newcommand{\cC}{\mathscr{C}}

\newcommand{\cI}{\mathscr{I}}

\newcommand{\calP}{\mathcal{P}}
\newcommand{\calS}{\mathcal{S}}

\newcommand{\calM}{\mathcal{M}}

\newcommand{\varp}{\varepsilon}

\newcommand{\rest}{{\restriction}}

\newenvironment{PROOF}[2][\proofname.]
   {\begin{proof}[#1]\renewcommand{\qedsymbol}{\textsquare$_{\mathrm{#2}}$}}
   {\end{proof}}

\renewcommand{\menos}{\smallsetminus}

\definecolor{bluet}{rgb}{0.0, 0.4, 0.6}

\hypersetup{
    colorlinks=true,
    linkcolor=bluet,
    filecolor=bluet,
    urlcolor=bluet,
    citecolor=bluet,
}

\definecolor{sub0}{RGB}{29,32,137}
\definecolor{sub1}{RGB}{1,71,157}
\definecolor{sub2}{RGB}{1,104,183}
\definecolor{sub3}{RGB}{0,160,234}
\definecolor{sug}{RGB}{0,154,68}
\definecolor{suy}{RGB}{208,219,1}

\definecolor{redun}{rgb}{0.65, 0.11, 0.19}
\definecolor{greenun}{rgb}{0.58, 0.71, 0.23}
\definecolor{dodger}{rgb}{0.0,0.5,1.0}
\definecolor{carrotorange}{rgb}{0.93, 0.57, 0.13}

\title[Absoluteness of the Riemann Integral]{Absoluteness of the Riemann Integral}

\author{Carlos M. Parra-Londoño}
\address{Universidad Nacional de Colombia, Campus El Volador, Departamento de Matemáticas, Facultad de Ciencias, Calle 59A $\#$ $63$-$20$, Medellín, Colombia}
\email{\href{mailto:cmparra@unal.edu.aco}{cmparra@unal.edu.co}}

\author{Andr\'es F.\ Uribe-Zapata}
\address{TU Wien, Faculty of Mathematics and Geoinformation, Institute of Discrete Mathematics and Geometry, Wiedner Hauptstrasse 8--10, A--1040 Vienna, Austria}
\email{\href{mailto:andres.zapata@tuwien.ac.at}{andres.zapata@tuwien.ac.at}}
\urladdr{\url{https://sites.google.com/view/andres-uribe-afuz}}

\thanks{Much of the content of this paper is derived from the second author's master's thesis. For this reason, the second author expresses profound gratitude to his primary advisor, Professor Diego A. Mejía from Kobe University, Japan, for his exceptional and meticulous guidance during the conception, development, and writing of the thesis. Additionally, the authors thank him for his review of this document, his assistance in refining its presentation, and for the invaluable suggestions that contributed to the results presented in this article. The second author was supported by the Austrian Science Fund (FWF): project number P33895.} 

\date{\today}



\begin{document}

\makeatletter
\def\@roman#1{\romannumeral #1}
\makeatother

\newcounter{enuAlph}
\renewcommand{\theenuAlph}{\Alph{enuAlph}}

\theoremstyle{plain}
  \newtheorem{theorem}{Theorem}[section]
  \newtheorem{corollary}[theorem]{Corollary}
  \newtheorem{lemma}[theorem]{Lemma}
  \newtheorem{mainlemma}[theorem]{Main Lemma}
  \newtheorem{maint}[theorem]{Main Theorem}
  \newtheorem{mainproblem}[theorem]{Main Problem}
  \newtheorem{construction}[theorem]{Construction}
  \newtheorem{prop}[theorem]{Proposition}
  \newtheorem{clm}[theorem]{Claim}
  \newtheorem{fact}[theorem]{Fact}
  \newtheorem{exercise}[theorem]{Exercise}
  \newtheorem{question}[theorem]{Question}
  \newtheorem{problem}[theorem]{Problem}
  \newtheorem{cruciallem}[theorem]{Crucial Lemma}
  \newtheorem{conjecture}[theorem]{Conjecture}
  \newtheorem{assumption}[theorem]{Assumption}
  \newtheorem{teorema}[enuAlph]{Theorem}

  \newtheorem*{corolario}{Corollary}
\theoremstyle{definition}
  \newtheorem{definition}[theorem]{Definition}
  \newtheorem{example}[theorem]{Example}
  \newtheorem{remark}[theorem]{Remark}
    \newtheorem{hremark}[theorem]{Historical Remark}
    \newtheorem{observation}[theorem]{Observation}
  \newtheorem{notation}[theorem]{Notation}
  \newtheorem{context}[theorem]{Context}

  \newtheorem*{defi}{Definition}
  \newtheorem*{acknowledgements}{Acknowledgements}

\numberwithin{equation}{theorem}
\renewcommand{\theequation}{\thetheorem.\arabic{equation}}

\def\sectionautorefname{Section}
\def\subsectionautorefname{Subsection}


\begin{abstract}

    This article explores the concept of \emph{absoluteness} in the context of mathematical analysis, focusing specifically on the Riemann integral on $\mathbb{R}^{n}$. In mathematical logic, \emph{absoluteness} refers to the invariance of the truth value of certain statements in different mathematical universes. Leveraging this idea, we investigate the conditions under which the Riemann integral on $\mathbb{R}^{n}$ remains absolute between transitive models of $\mathrm{ZFC}$ —the standard axiomatic system in which current mathematics is usually formalized. To this end, we develop a framework for integration on Boolean algebras with respect to finitely additive measures and show that the classical Riemann integral is a particular case of this generalized approach. Our main result establishes that the Riemann integral over rectangles in $\mathbb{R}^{n}$ is absolute in the following sense: if $M \subseteq N$ are transitive models of $\mathrm{ZFC}$, $a, b \in \mathbb{R}^{n} \cap M$, and $f \colon [a, b] \to \mathbb{R}$ is a bounded function in $M$, then $f$ is Riemann integrable in $M$ if, and only if, in $N$ there exists some Riemann integrable function $g \colon [a, b] \to \mathbb{R}$ extending $f$. In this case, the values of the integrals computed in each model are the same. Furthermore, the function $g$ is unique except for a measure zero set.
\end{abstract}

\maketitle

{
\small
\hypersetup{linkcolor=black}
\tableofcontents
}

\section{Introduction}

In mathematical logic, the concept of \emph{absoluteness} describes the invariance of the validity of certain statements across different universes of mathematics —or fragments of mathematics— that we call \emph{models}. Broadly speaking, a statement is \emph{absolute} relative to two models if it retains its validity in both, reflecting a logical stability that does not depend on the extensions of the model (see~\autoref{n.2}). For example, the so-called \emph{$\Delta_{0}$ formulas} —those whose quantifiers are bounded— are inherently absolute. Consequently, many fundamental notions in set theory, such as the empty set, union, intersection, Cartesian product of sets, being a function, an ordinal, and even a natural number, are absolute between any transitive model and $V$, the universe of set theory. Although $\Delta_{0}$ formulas are too elementary in a certain sense, there are results that guarantee the absoluteness of more complex formulas, such as the \emph{Mostowski Absoluteness Theorem} and the \emph{Shoenfield's Absoluteness Theorem} (see e.g.~\cite{MejRiv}). The first one states that any analytic subset of the Baire space $\bbN^{\bbN}$ is absolute for transitive models of $\ZFC$ —the standard axiomatic system in which current mathematics is usually formalized. The second one states that any $\Sigma_{2}^{1}$ subset of the Baire space is absolute for certain models of $\ZFC$ (see~\cite{Kechris}).  Beyond these results, the absoluteness of more complex formulas is often related to the existence of so-called \emph{large cardinals} (see e.g.~\cite{BagariaF} and~\cite{Bagaria2}).

The notion of absoluteness transcends the boundaries of set theory and model theory, significantly influencing areas such as number theory, geometry, topology, and even the philosophy of mathematics. For instance, absolute definitions of natural numbers ensure their stability across a wide variety of models, providing a robust foundation for arithmetic analysis even in more general contexts, such as non-standard analysis, which has applications in many diverse areas of mathematics (see e.g~\cite{Loeb2000}). In geometry and topology, the notion of absoluteness enables the identification of certain geometric properties that remain invariant with respect to the considered universe. This makes the concept a valuable tool in the interaction between geometry and model theory (see, e.g.~\cite{Chang1990} and~\cite{MacLane1992}). From a philosophical perspective, the notion of absoluteness addresses fundamental questions about the nature of mathematical truth, allowing a distinction between properties that depend on particular axioms and those that have a universal and invariant character. This approach provides essential tools for examining the relationship between axioms, definability, and provability (see e.g.~\cite{Lemanska2015},~\cite{Krol2021}, and~\cite{Halbach2014}).

In this paper, we focus on studying the absoluteness of a concept from mathematical analysis: the Riemann integral. This is unexpectedly motivated by problems related to forcing theory and consistency proofs. Specifically, in 2000, Saharon Shelah proved in~\cite{Sh00} that a certain cardinal invariant, called the \emph{covering number of the null ideal}, may have countable cofinality, thus solving a problem posed by David Fremlin that had remained open for nearly 30 years (see~\cite[Sec.~5.1]{uribethesis}). To achieve this, Shelah introduced a new forcing technique which —without being too technical— involves the use of finite-support iterations accompanied by finitely additive measures on $\calP(\bbN)$, which satisfy certain special conditions. Recently, based on Shelah's aforementioned work and another article by Jakob Kellner, Saharon Shelah, and Anda T\u{a}nasie (see~\cite{KST}), Miguel A. Cardona, Diego A. Mejía, and the second author, in~\cite{CMU}, introduced a general theory of iterated forcing using finitely additive measures\footnote{A preliminary version of this general theory of iterated forcing with finitely additive measures was introduced in the master's thesis of the second author (see~\cite{uribethesis}). In this thesis, an entire chapter was dedicated to the study of integration on Boolean algebras, and some results concerning the absoluteness of this integral were proven, upon which this article is based (see \cite[Ch.~3]{uribethesis}).}. To prove one of their main results —specifically, the theorem of extension at successor steps (see~\cite[Sec.~7.2]{CMU})— it was necessary to rely on the absoluteness of the integral for certain functions in ${}^{\bbN} \bbN$ with respect to  finitely additive measures on $\calP(\bbN)$.  This result was not difficult to obtain, as the functions in question were defined on $\mathbb{N}$, which is absolute for transitive models of $\mathsf{ZFC}$ (see~\autoref{n.2}). However, the situation is quite different in the case of the Riemann integral. Since the set of real numbers is not absolute, if we have two transitive models $M \subseteq N$ of $\mathsf{ZFC}$, and in $M$ we have a bounded real-valued function defined on some rectangle $[a, b]$, then interpreting $f$ in $N$ may result in $f$ not being defined on the whole  $[a, b]$, as new real numbers may appear in $N$. This makes the analysis of the absoluteness of the Riemann integral a more complex problem. In this article, we address it and establish the following main result, which corresponds to~\autoref{s20}.

\begin{teorema}\label{ta}
    Let $M, N$ be transitive models of $\ZFC$ such that $M \subseteq N$, and $n \in \bbN$. Let $a = (a_{0}, \dots, a_{n-1}) $ and $ b = (b_{0} \dots, b_{n-1})$ be in $ \bbR^{n} \cap M$ such that for any $i < n$, $a_{i} \leq b_{i}$, and $[a, b] \coloneqq \prod_{i < n} [a_{i}, b_{i}]$.  In $M$, assume that $f \colon [a, b] \to \bbR$ is a bounded function. Then, $f$ is Riemann integrable in $M$ if, and only if, in $N$ there exists some Riemann integrable function $g$ on $[a, b]$ extending $f$.  In this case, 
     $$
    \left( \int_{[a, b]} f   \right)^{M} = \left(  \int_{[a, b]} g  \right)^{N} \! \! \! \!, 
    $$ where the left value is the result of the Riemann integral computed in $M$, and the right one is the result of the integral computed in $N$.  Moreover, $g$ is unique except in a Lebesgue measure zero set: if in $N$, $g^{\ast}$ is another Riemann integrable function on $[a, b]$ extending $f$, then there exists some measure zero set $E \subseteq [a, b]$ such that for any $x \in [a, b] \setminus E$, $g^{\ast}(x) = g(x)$. 
\end{teorema}

To prove \autoref{ta}, in the first part of~\autoref{3.1}, we review some elementary properties of finitely additive measures on Boolean algebras, and in the second part, we introduce  the needed basic notions of absoluteness. Finally, in~\autoref{3}, following the ideas of the Riemann integral, and based on~\cite[Ch.~3]{uribethesis} and~\cite{CMUP} —an ongoing project by Miguel A. Cardona, Diego A. Mejía, and the second author, which conducts an in-depth study of finitely additive measures and their integrals on Boolean algebras— we define a notion of integral on Boolean algebras with respect to finitely additive measures, we prove that this integral is absolute for transitive models of $\ZFC$ (see~\autoref{t100}), and that the Riemann integral over rectangles in $\bbR^{n}$ is a particular case of this integral (see~\autoref{r20}). This framework will allow us to prove several results related to the absoluteness of our integral over Boolean algebras and, in particular, prove our main result~\autoref{ta}. 

\section{Preliminaries}\label{3.1}

In this section, we will introduce the basic notions of finitely additive measures on Boolean algebras, study some of their fundamental properties, and examine their relationship with filters and ultrafilters.

We begin by introducing some notation and essential concepts of Boolean algebras.

\subsection{Boolean algebras}

Recall that $\cB \coloneqq \langle \cB, \wedge, \vee, \neg, 0_{\cB}, 1_{\cB} \rangle$ is a \emph{Boolean algebra} if $\cB$ is a non-empty set, $\wedge, \vee$ are binary commutative and associative operations on $\cB,$ $\neg$ is a unary operation on $\cB,$ $1_{\cB}, 0_{\cB} \in \cB$, and the following properties are satisfied for all $a, b, c \in \cB$: 

    \begin{multicols}{2}
        \begin{enumerate}    
            \item \emph{Absorption}: 
            \index{Boolean algebra!absorption}
    
            \begin{enumerate}[label = \normalfont (\roman*)]
                \item $(a \vee b) \wedge b = b$,
    
                \item $(a \wedge b) \vee b = b.$
            \end{enumerate}
    
            \item \emph{Distributivity}:             \index{Boolean algebra!distributivity}
    
            \begin{enumerate}[label = \normalfont (\roman*)]
                \item $(a \vee b) \wedge c = (a \wedge c) \vee (b \wedge c),$
    
                \item $(a \wedge b) \vee c = (a \vee c) \wedge (b \vee c).$
            \end{enumerate}
    
            \item \emph{Identity}:       \index{Boolean algebra!identity}
    
            \begin{enumerate}[label = \normalfont (\roman*)]
                \item $a \wedge 1_{\cB}  = a,$
    
                \item $a \vee 0_{\cB} = a.$  
            \end{enumerate}
    
            \item \emph{Complementation}:  \index{Boolean algebra!complements}
    
            \begin{enumerate}[label = \normalfont (\roman*)]
                \item $a \vee \neg  a = 1_{\cB},$
    
                \item $a \wedge \neg a = 0_{\cB}.$ 
            \end{enumerate}
        \end{enumerate}
    \end{multicols}

The operations $\wedge$ and $\vee$ are known as \emph{meet} and \emph{join} respectively. Also $\wedge, \vee$ and $\neg$ are known as the \emph{Boolean operations of $\cB$}. 

The canonical example of a Boolean algebra is the power set: Consider a non-empty set $X.$ Then $\langle \calP(X),  \cap,  \cup,  {}^{\rm{c}},  \emptyset,  X \rangle$ is a Boolean algebra, where ${}^{\mathrm{c}} \colon \calP(X) \to \calP(X)$ is defined by $A^{\mathrm{c}} \coloneqq X \setminus A$. 

\begin{definition}\label{s14}
    Let $\cB$ be a Boolean algebra and $a, b \in \cB$. 
    \begin{multicols}{2}
        \begin{enumerate}
            \item $a \menos b \coloneqq a \wedge \neg b$. 
    
            \item $a \leq_{\cB} b $ iff $a \wedge b = a.$
    
            \item $\cB^{+} \coloneqq \cB \setminus \{ 0_{\cB} \}$. 

            \item $a$, $b$ are \emph{incompatible} iff $a \wedge b = 0_{\cB}$.
        \end{enumerate}
    \end{multicols}
\end{definition}
 
Notice that, $\menos$ is a binary operation on $\cB$ and $\leq_{\cB}$ is a partial order on $\cB$, which allows us to characterize the values of the operations $\vee$ and $\wedge$: if $a, b \in \cB$, then $a \wedge b$ and $a \vee b$ are the \emph{infimum} and the \emph{supremum} of $\{ a, b \}$, respectively, in the sense that they are the minimum upper bound and the maximum lower bound, respectively, with respect to the order $\leq_{\cB}$. If $I$ is an arbitrary set, and $\{ b_{i} \colon i \in I \} \subseteq \cB$, we can define $\bigvee_{i \in I} b_{i}$ and $\bigwedge_{i \in I} b_{i}$ in a similar way but, in general, existence is only guaranteed when $I$ is finite. A Boolean algebra such that every infinite subset has a supremum is called \emph{complete}, and a Boolean algebra in which every countable non-empty subset has supremum is called \emph{$\sigma$-complete}.

Let $\cB$ be a Boolean algebra. A \emph{Boolean sub-algebra} of $\cB$ is a non-empty subset $\cC \subseteq \cB$ that is closed under the Boolean operations of $\cB$. Consequently, $\cC$ is itself a Boolean algebra, and it contains both $0_{\cB}$ and $1_{\cB}$. 

Recall that a \emph{Boolean homomorphism} is a function $h \colon \cB \to \cC$ —where $\cB$ and $\cC$ are Boolean algebras— that preserves the Boolean operations, that is,  $h(a \wedge b) = h(a) \wedge h(b),$ $h(a \vee b) = h(a) \vee h(b)$, and   $h(\neg a) =  \neg h(a)$ for all $a, b \in \cB$, which implies that $h(0_{\cB}) = 0_{\cC}$ and $h(1_{\cB}) = 1_{\cC}$. A \emph{Boolean isomorphism} from $\cB$ into $\cC$ is a bijective Boolean homomorphism from $\cB$ onto $\cC.$ 

\begin{example}\label{s16}
    Let $X$, $Y$ be non-empty sets and $h \colon X \to Y$ a function. 
    
    \begin{enumerate}
        \item The map $f_{h} \colon \calP(Y) \to \calP(X)$ defined by $f_{h}(A) \coloneqq h^{-1}[A]$ for every $A \subseteq Y$ is a homomorphism. Furthermore, $f_{h}$ is an isomorphism if, and only if, $h$ is a bijection. 

        \item More generally, if $\cB$ is a Boolean sub-algebra of $\calP(X)$, then $\cC \coloneqq h^{\to}(\cB)$ is a Boolean sub-algebra of $\calP(Y)$, where $ h^{\to}(\cB) \coloneqq \{ A \subseteq Y \colon h^{-1}[A] \in \cB \}$, and the map $f_{h} \colon \cC \to \cB$ defined by $f_{h}(c) \coloneqq h^{-1}[c]$ for every $c \in \cC$ is a homomorphism. Furthermore, $f_{h}$ is one-to-one iff $h$ is onto, and if $F[\cB] = \cC$, then $h$ is onto, where $F \colon \calP(X) \to \calP(Y)$ is defined by $F(A) \coloneqq h[A]$ for every $A \subseteq X$. 
    \end{enumerate}
\end{example}

As a consequence of \emph{Stone's representation theorem} (see \cite[Thm.~4.1]{BellM}), we can characterize any Boolean algebra using $\calP(X)$ for some set $X$. 

\begin{theorem}\label{b40}
    Every boolean algebra is isomporhic to a Boolean sub-algebra of $\calP(X)$ for some set $X.$  
\end{theorem}

The core of the proof of Stone's representation theorem is the notion of \emph{ultrafilter}, which we introduce in the following definition.  

\begin{definition}
    Let $\cB$ be a Boolean algebra.
    
    \begin{enumerate}
        \item A \emph{filter on $\cB$} is a non-empty set $F \subseteq \cB$ such that:  

            \begin{enumerate}[label = \normalfont (\roman*)]
                \item if $x, y \in F,$ then $x \wedge y \in F,$
            
                \item if $x \in F$ and $x\leq_{\cB} y,$ then $y \in F,$
            
                \item $0_{\cB} \notin F.$  
            \end{enumerate}

        \item An \emph{ultrafilter on $\cB$} is a filter $F \subseteq \cB$ such that, for any $b \in \cB, $ either $b \in F$ or $\neg b \in F.$ 

        \item We say that a non-empty set $I \subseteq \cB$ is an \emph{ideal on $\cB$} if it satisfies the following conditions: 

            \begin{enumerate}[label = \normalfont (\roman*)] 
                \item if $x, y  \in I,$ then $x \vee  y \in I,$ 
            
                \item if $x \in I$ and $y \leq_{\cB} x,$ then $y \in I,$ 
            
                \item $1_{\cB} \notin I.$
            \end{enumerate} 
    \end{enumerate}
\end{definition}

Filters and ideals on a Boolean algebra are \emph{dual notions} in the following sense. 

\begin{fact}\label{b65}
    Let $\cB$ be a Boolean algebra and $F \subseteq \cB$. Define  $F^{\neg} \coloneqq \{ \neg a\colon a \in F \}.$  Then $F$ is a filter on $\cB$ if, and only if, $F^{\neg}$ is an ideal on $\cB.$ \index{$F^{\neg}$}
\end{fact}

It is straightforward to show that the intersection of all Boolean sub-algebras of $\cB$ containing a given subset $B$ is itself a Boolean sub-algebra of $\cB$. This sub-algebra, known as the \emph{Boolean sub-algebra generated by $B$}, is the smallest —with respect to $\subseteq$— Boolean sub-algebra of $\cB$ containing $B$. This is denoted by $\langle B \rangle_{\cB}$, or simply $\langle B \rangle$ when the context is clear. In this setting, $B$ is referred to as the \emph{generating set} of $\langle B \rangle_{\cB}$. 

Since filters are upwards closed and closed under $\wedge,$ and ideals are downwards closed and closed under $\vee,$ we can characterize the Boolean sub-algebra generated by a filter as follows.

\begin{fact}\label{b54}
    Let $\cB$ be a Boolean algebra. If $F \subseteq \cB$ is a filter, then $\langle F \rangle = F \cup F^{\neg}.$ As a consequence, $F$ is an ultrafilter on $\cB$ iff $\langle F \rangle = \cB.$ 
\end{fact}

We close this subsection by introducing the notion of \emph{free filter}. 

\begin{definition}\label{b67}
    Let $X$ be a non-empty set. We say that $F \subseteq \calP(X)$ is a \emph{free filter} if it is a filter containing all the co-finite subsets of $X$. Equivalently, if $F^{\neg}$ is an ideal including all finite subsets of $X$. 
\end{definition}

\subsection{Finitely additive measures on Boolean algebras}

Without resorting to the notion of $\sigma$-algebra, we can generalize the idea of \emph{measure} in the context of Boolean algebras: a measure on a Boolean algebra $\cB$ is a function $\bfm \colon \cB \to [0, \infty]$ such that $\bfm(0_{\cB}) = 0$ and, if $\{ b_{n} \colon n \in \bbN\} \subseteq \cB$ is such that $\bigvee_{n \in \bbN} b_{n} \in \cB,$ then $$ \bfm \left( \bigvee_{n \in \bbN} b_{n} \right) = \sum_{n \in \bbN} \bfm(b_{n}), $$ whenever for any $i, j \in \bbN,$ if $i \neq j,$ then $b_{i} \wedge b_{j} = 0_{\cB}.$  If we weaken this condition and enforce it only for finite sets of elements, we get the notion of \emph{finitely additive measures} on Boolean algebras. 

\begin{definition}\label{m4}    
    Let $\cB$ be a Boolean algebra. A \emph{finitely additive measure on $\cB$} is a function $\Xi \colon \cB \to [0,\infty]$ satisfying:
    
    \begin{enumerate}[label = \normalfont (\roman*)]
        \item $\Xi (0_{\cB})=0$,
        
        \item $\Xi(a\vee b)=\Xi(a)+\Xi(b)$ whenever $a,b\in\cB$ and $a \wedge b= 0_{\cB}$.
    \end{enumerate}

    We say that $b \in \cB$ has \emph{$\Xi$-measure} $\delta$ if $\Xi(b) = \delta.$ In general, we exclude the \emph{trivial finitely additive measure}, that is,  when talking about finitely additive measures, we will always assume $\Xi(1_\cB) > 0.$ We will occasionally use the acronym ``fam'' or ``FAM'' to refer to finitely additive measures. 
\end{definition}

There are several types of finitely additive measures. In the following definition, we introduce some of them. 

\begin{definition}\label{m9}
    Let $\cB$ be a Boolean algebra and $\Xi$ a finitely additive measure on $\cB.$ Then:

    \begin{enumerate}[label = \normalfont (\arabic*)]
        \item   We say that $\Xi$ is \emph{finite} if $\Xi(1_{\cB}) < \infty.$
    
        \item When $\Xi(1_\cB) = 1$, we say that $\Xi$ is \emph{a probability finitely additive measure}.
        
        \item If $\Xi(b)>0$ for any $b \in \cB^{+},$ we say that $\Xi$ is \emph{strictly positive}.

        \item\label{m9.4} If $\cB$ is a Boolean sub-algebra of $\mathcal{P}(X)$, we say that $\Xi$ is a \emph{free} if, for any $x\in X$, $\{x\}\in\cB$ and $\Xi(\{x\}) = 0$.
    \end{enumerate}
\end{definition}

We adopt the name \emph{free finitely additive measure} in connection with the notion of \emph{free filter} (see~\autoref{b67}).

\begin{example}\label{m37}\ 

    \begin{enumerate}
        \item Let $\cB$ be a Boolean algebra, $\Xi$ a finitely additive measure on it, and $b \in \cB$ with positive finite measure. We define the function 
        $\Xi_{b} \colon \cB \to [0,1]$ by $\Xi_{b}(a) \coloneqq \frac{\Xi(a \wedge b)}{\Xi(b)}$ for any $a \in \cB$. It is clear that $\Xi_{b}$ is a finitely additive probability measure. 
        
        \item\label{m37.b} Let $X$ be a non-empty set. For a finite non-empty set $u \in \calP(X)$, we define $\Xi^u \colon \calP(X) \to [0, 1]$ by $\Xi^{u}(x) \coloneqq \frac{\vert x \cap u \vert}{ \vert u \vert}$ for any $x \in \calP(X)$. We call this the \emph{uniform measure with support $u$}. 
    \end{enumerate}    
\end{example}

In general $\Xi^{u}$ is not strictly positive. To guarantee the existence of more interesting finitely additive measures, we must require that the Boolean algebra has more structure, for instance, to be $\sigma$-centered. Recall that a Boolean algebra $\cB$ is $\sigma$-centered whenever $\cB^{+}$ can be decomposed as a countable union of ultrafilters on $\cB$. 

\begin{theorem}\label{m38}
    Every $\sigma$-centered Boolean algebra admits a strictly positive probability finitely additive measure.  
\end{theorem}

\begin{PROOF}[\textbf{Proof}]{\ref{m38}}
    Let $\cB$ be a $\sigma$-centered Boolean algebra witnessed by $\{ F_{n} \colon n \in \bbN\}$. For any $b \in \cB,$ consider $\omega_{b} \coloneqq \{ n \in \bbN\colon b \in  F_{n} \}$ and set $\Xi \colon  \cB \to [0, 1]$ such that 
    $$
    \Xi(b) \coloneqq  \sum_{n \in \omega_{b}} \frac{1}{2^{n+1}}. 
    $$
    
    Let $a, b \in \cB^{+}$ such that $a \wedge b = 0_{\cB}.$ It is clear that $\omega_{a} \cap \omega_{b} = \emptyset$ and $\omega_{a} \cup \omega_{b} \subseteq \omega_{a \vee b}.$ Conversely, let $m \in \omega_{a \vee b},$ so $a \vee b \in F_{m}.$ If $m \notin \omega_{a}$ and $m \notin \omega_{b},$ then $\neg a\in F_{m}$ and $\neg  b \in F_{m},$ and therefore, $\neg  (a \vee b) = \neg a \wedge \neg b  \in F_{m},$ which is a contradiction. Thus, $\omega_{a} \cup \omega_{b} = \omega_{a \vee b}$, and we can calculate 
    $$
    \Xi(a \vee b) = \sum_{n \in \omega_{a \vee b}} \frac{1}{2^{n+1}} = \sum_{n \in \omega_{a}} \frac{1}{2^{n+1}} + \sum_{n \in \omega_{b}} \frac{1}{2^{n+1}} = \Xi(a) + \Xi(b).
    $$ 
    
    Finally, it is clear that $\Xi$ is strictly positive and, since $1_{\cB} \in F_{n}$ for all $n \in \bbN, \omega_{1_{\cB}} = \bbN,$ hence $\Xi(1_{\cB}) = 1,$ that is, $\Xi$ is a  finitely additive probability measure. 
\end{PROOF}

In general, proving the existence of interesting finitely additive measures requires the Axiom of Choice, which relies on non-constructive methods (see \cite{Lauwers}). In the next section, we will explore the close relationship between finitely additive $\{0, 1\}$-valued measures and ultrafilters, from which finitely additive measures naturally arise as examples (see~\autoref{m70}).

Next, we show some elementary properties of finitely additive measures.

\begin{lemma}\label{n2}
    Let $\cB$ a Boolean algebra,  $\Xi$ a finitely additive measure on $\cB$, and $a, b \in \cB$. Then: 
    \begin{enumerate}[label = \normalfont (\alph*)]
        \item\label{n2.a}  If $a \leq_{\cB} b$ then $\Xi(a) \leq \Xi(b).$ 

        \item\label{n2.b} $\Xi(a \vee b) + \Xi(a \wedge b) = \Xi(a) + \Xi(b)$.

        \item\label{n2.c} If $\Xi(a \vee b) < \infty$, then $\Xi(a) - \Xi(b) \leq \Xi(a \setminus b)$ and $\Xi(b) - \Xi(a) \leq \Xi(b \setminus a)$. 

        \item\label{n2.d} If $n \in \bbN$ and $\langle b_{i} \colon i < n \rangle \subseteq \cB$, then 
        $$\Xi \left ( \bigvee_{i < n} b_{i} \right) \leq \sum_{i < n} \Xi(b_{i}).$$
        Furthermore, the equality holds whenever $\langle b_{i} \colon i < n  \rangle$ is a sequence of pairwise incompatible elements of $\cB$.  

        \item\label{n2.e} $\Xi(1_{\cB}) = \Xi(b) + \Xi(\neg  b)$.
    \end{enumerate}
\end{lemma}

\begin{PROOF}{\ref{n2}}
    \ref{n2.a}: Assume that $a \leq_{\cB} b.$ Since $b = a \vee (b \menos a)$ and $a \wedge (b \menos a) = 0_{\cB},$ we have that $\Xi(a) \leq \Xi(a) + \Xi(b \menos a) = \Xi(b).$ Thus $\Xi(a) \leq \Xi(b).$  

    \ref{n2.b}: Since $a = (a \wedge b) \vee (a \menos b)$ and $(a \wedge b)  \wedge (a \menos b) = 0_{\cB},$ we get $\Xi(a) = \Xi(a \wedge b) + \Xi(a \menos b)$. Similarly,  we have that $\Xi(b) = \Xi(b \menos a) + \Xi(b \wedge a).$ As a consequence, 
    \begin{equation*}
        \begin{split}
            \Xi(a) + \Xi(b) & = [\Xi(a \wedge b) + \Xi(a \menos b) + \Xi(b \menos a)] + \Xi(a \wedge b) = \Xi(a \vee b) + \Xi(a \wedge b). 
        \end{split}
    \end{equation*}

    \ref{n2.c}: Follows from~\ref{n2.b} since $\Xi(a) - \Xi(b) = \Xi(a \menos b) - \Xi(b \menos a)$.

    \ref{n2.d}: By an inductive argument, it is enough to prove it only for two elements $b_{0}, b_{1} \in \cB$. By~\ref{n2.b}, $\Xi(b_{0} \vee b_{1}) \leq \Xi(b_{0} \vee b_{1}) + \Xi(b_{0} \wedge b_{1}) =  \Xi (b_{0}) + \Xi(b_{1})$ and the equality holds when $b_{0} \wedge b_{1} = 0_{\cB}.$  

    \ref{n2.e}: Straightforward from~\ref{n2.b}. 
\end{PROOF}




We close this section by showing that we can transfer finitely additive measures using Boolean homomorphisms. 

\begin{definition}\label{ss0}
    \
    \begin{enumerate}
        \item Let $\cB$, $\cC$ be Boolean algebras and $f \colon \cC \to \cB$ a Boolean homomorphism. If $\Xi$ is a finitely additive measure on $\cB,$  define the \emph{finitely additive measure induced by $f$ and $\Xi$ on $\cC$}, denoted $\Xi_{f}$, by $\Xi_{f
        }(c) \coloneqq \Xi(f(c))$ for all $c \in \cC$.  

        \item Let $X$, $Y$ be non-empty sets, $h \colon X \to Y$, $\cB$ a Boolean sub-algebra of $\calP(X)$, $\Xi$ a finitely additive measures on $\cB$, and $\cC \coloneqq h^{\to}(\cB)$. We define $\Xi_{h} \coloneqq \Xi_{f_{h}}$, where $f_{h}$ is as in~\autoref{s16}.  
    \end{enumerate}
\end{definition}

\begin{fact}
    $\Xi_{h}$ is a finitely additive measure. Moreover, $\Xi_{h}$ is a probability iff $\Xi$ is a probability. Furthermore, if for any $x \in X$, $\{ x \} \in \cB$, and $h$ is finite-to-one, then $\Xi_{h}$ is free iff $\Xi$ is free. 
\end{fact}

\subsection{Connections with filters and ultrafilters}\label{3.2}

In this section, we study the connection that exists between the $\{0, 1\}$-valued finitely additive measures and the ultrafilters on a Boolean algebra. 

We start by showing that every filter $F$ naturally induces a finitely additive probability measure on $\langle F \rangle$. 

\begin{lemma}\label{m70}
    Let $\cB$ be a Boolean algebra and $F$  a filter on $\cB.$ Then $\Xi_{F} \colon \langle F \rangle \to \{0, 1\}$ such that, for any $b \in \langle F \rangle,$ 
    $$\Xi_{F}(b)= \left\{ \begin{array}{lcc}
             1 &   \rm{if}  & b \in F,  \\
             
             \\ 0 &  \rm{if}  & b \in F^{\neg}, 
             \end{array}
   \right.$$ 
   
   is a finitely additive probability measure. Furthermore, if $G$ is another filter on $\cB,$ then  $$F \subseteq G \Leftrightarrow \Xi_{F} \leq \Xi_{G}.$$ Moreover, if $F \subseteq \calP(X)$ for some set $X$ then $F$ is free if, and only if, $\Xi_{F}$ is free. 
\end{lemma}

\begin{PROOF}{\ref{m70}}
    Notice that $\Xi_{F}$ is well-defined because, by \autoref{b54}, $\langle F \rangle = F \cup F^{\neg}$ and those sets are disjoint. Since $F$ is a filter, by \autoref{b65}, $F^{\neg}$ is an ideal, so $0_{\cB} \in F^{\neg},$ hence $\Xi_{F}(0_{\cB}) = 0.$ To show that $\Xi_{F}$ is a finitely additive measure, let $a, b \in \langle F \rangle$ such that $a \wedge b = 0_{\cB}$. If $a \in F, $ then $b \in F^{\neg},$ hence $\Xi_{F}(a \vee b) = 1  = \Xi_{F}(a) + \Xi_{F} (b).$ The case $a \in F^{\neg}$ and $b \in F$ is analogous. If $a, b \in F^{\neg},$ then $a \vee b \in F^{\neg},$ that is, $\Xi_{F}(a \vee b) = 0 = \Xi_{F}(a) + \Xi_{F}(b).$  Thus, $\Xi_{F}$ is a finitely additive measure on $\langle F \rangle$, and clearly $\Xi_{F}(1_{\cB}) = 1$.  

    Now, assume that $F \subseteq G$ and let  $b \in \langle F \rangle.$ On the one hand,  if $b \in F,$ then $b \in G$ and therefore, $\Xi_{F}(b) = 1 = \Xi_{G}(b).$ On the other hand, if $b \notin F,$ then $\Xi_{F}(b) = 0 \leq \Xi_{G}(b).$ Thus, in any case $\Xi_{F}(b) \leq \Xi_{G}(b).$ Conversely, assume that $\Xi_{F} \leq \Xi_{G}$ and let $b \in F,$ hence $1 = \Xi_{F}(b) \leq \Xi_{G}(b),$ therefore $\Xi_{G}(b) = 1,$ that is, $b \in G.$  Thus, $F \subseteq G.$ 

    Finally, that $F$ is free if and only if $\Xi_{F}$ is free follows directly from the definitions of freeness (see~\autoref{b67} and~\autoref{m9}~\ref{m9.4}).
\end{PROOF}

If we choose a suitable ultrafilter, we can use~\autoref{m70} to construct an example of a finitely additive measure that is not a measure.

\begin{example}\label{m72}
    Let $\cB$ be a Boolean sub-algebra of $\calP(X),$ where $X$ is a countable set, and let $F \subseteq \cB$ be an ultrafilter on $\cB.$ If $F$ is a free on $\cB,$ then $\Xi_{F}$ is a finitely additive measure on $\cB$ that is not a measure on $\cB.$ Indeed, assume that $F$ is a free on $\cB.$ We already know that $\Xi_{F}$ is a finitely additive measure on $\cB$ by \autoref{b54} and \autoref{m70}. Now, since $X$ is countable, we can write $X = \{ x_{n} \colon n \in \bbN\}.$ For any $n \in \bbN,$ define $B_{n} \coloneqq \{ x_{n} \}.$ Notice that $B_{n}^{\rm{c}}$ is co-finite, and therefore, $B_{n}^{\rm{c}} \in F,$ that is, $B_{n} \in F^{\neg}.$ Thereby, $\Xi_{F}(B_{n}) = 0$ for any $n \in \bbN.$ However, $X \in F$ because it is co-finite, hence $\Xi_{F}(X) = 1.$ Thus, $$\Xi_{F} \left( \bigcup_{n \in \bbN} B_{n} \right) = \Xi_{F}(X) = 1 \neq 0 = \sum_{n \in \bbN} \Xi_{F}(B_{n}).$$ Consequently, $\Xi_{F}$ is not a measure on $\cB$.   
\end{example}

Conversely (see~\autoref{m70}), finitely additive probability measures also induce filters and in some cases, ultrafilters. 

\begin{lemma}\label{m74}    
    Let $\cB$, $\cC$ be Boolean algebras with $\cC \subseteq \cB$, and $\Xi \colon \cC \to \{0, 1 \}$ a finitely additive probability measure. Define $F_{\Xi} \coloneqq \{ c \in \cC \colon \Xi(c) = 1  \}$. Then:

    \begin{enumerate}[label = \normalfont (\alph*)]
        \item\label{m74.a} $F_{\Xi}$ is a filter on $\cC$.

        \item\label{m74.g} If $\rho$ is a finitely additive probability measure on $\cC$ then $F_{\Xi} \subseteq F_{\rho}$ iff $\Xi \leq \rho$. 

        \item\label{m74.b} $ \langle F_{\Xi} \rangle = \cC$.

        \item\label{m74.c} $\Xi_{F_{\Xi}} = \Xi,$ where $\Xi_{F_{\Xi}}$ is as in~\autoref{m70}.

        \item\label{m74.e} $F_{\Xi_{F}} = F$, where $\Xi_{F}$ is as in~\autoref{m70}.

        \item\label{m74.d} $F_{\Xi}$ is an ultrafilter on $\cB$ iff $\cB = \cC$.  
    \end{enumerate}
\end{lemma}

\begin{PROOF}{\ref{m74}}
    \ref{m74.a}: Let $c, d \in F_{\Xi}$. By~\autoref{n2}~\ref{n2.b}, $ \Xi(c \wedge d) =  2 - \Xi(c \vee d)$, and therefore $\Xi(c \wedge d) = 1.$ Thus, $c \wedge d \in F_{\Xi}.$ Now, if $a \in \cC$ and $c \leq_{\cC} a$, by~\autoref{n2}~\ref{n2.a}, $\Xi(a) \geq \Xi(c) = 1,$ hence $\Xi(a) = 1.$ Thus, $a \in F_{\Xi}.$ Finally, since $\Xi(0_{\cB}) = 0,$ it follows that $0_{\cB} \notin F_{\Xi}.$ Thus, $F_{\Xi}$ is a filter on $\cC.$ 

    \ref{m74.g}: Let $c \in \cC$. If $\rho(c) < \Xi(c)$ then $c \in F_{\Xi}$ and $c \notin F_{\rho}$, which shows that $F_{\Xi} \subseteq F_{\rho}$ implies $\Xi \leq \rho$. The converse is clear. 

    \ref{m74.b}: By \autoref{n2}~\ref{n2.d}, if $c \in \cC$ and $\Xi(c) = 0$ then $\Xi(\neg  c) = \Xi(1_{\cC}) - \Xi(c) = 1,$ hence $\neg  c \in F_{\Xi},$ that is, $c \in F_{\Xi}^{\neg}.$ Thus, $\cC = F_{\Xi} \cup F_{\Xi}^{\neg} = \langle F_{\Xi} \rangle.$ 

    \ref{m74.c}: Since $\cC = \langle F_{\Xi} \rangle, \, \Xi$ and $\Xi_{F_{\Xi}}$ have the same domain. For $c \in \cC,$ we have that  $\Xi_{F_{\Xi}}(c) = 1 \Leftrightarrow c \in F_{\Xi} \Leftrightarrow \Xi(c) = 1. $ Thus, $\Xi_{F_{\Xi}} = \Xi.$ 

    \ref{m74.e}: It is clear since, for any $c \in \cC$, $c \in F$ iff $\Xi_{F}(c) = 1$ iff $c \in F_{\Xi_{F}}$. 

    \ref{m74.d}: Straightforward. 
\end{PROOF}

In general $F_{\Xi}$ is not a filter on $\cB$. To obtain this, we need to close it upwards.

\begin{corollary}\label{m78}
    Let $\cB$ be a Boolean algebra and $\cC$ a Boolean sub-algebra of $\cB.$ Then every finitely additive probability measure $\Xi \colon \cC \to \{0, 1\}$ induces a filter on $\cB,$ namely, $F_{\Xi}^{\uparrow} \coloneqq \{ b \in \cB \colon \exists c \in F_{\Xi}(c \leq_{\cB} b) \}$.
\end{corollary}

As a consequence of~\autoref{m74} and~\autoref{m70}, we can establish some relations between the collections of filters, ultrafilters and $\{ 0,1 \}$-valued finitely additive probability measures associated to a Boolean algebra. However, we need to introduce some notation first.  

\begin{definition}
    For a given Boolean algebra $\cB$, let $\mathcal{UF}_{\!\cB}$ be the collection of all ultrafilters on $\cB$. Notice that $\mathcal{UF}_{\! \cB}$ is partially ordered by inclusion. Similarly, let $\mathcal{FAM}_{0, 1, \cB}$ be the collection of all $\{ 0, 1 \}$-valued finitely additive probability measures on $\cB$. Notice that, $\mathcal{FAM}_{0, 1, \cB}$ is partially ordered by the point-wise order on functions $\leq$.
\end{definition}

If we apply~\autoref{m74} to ultrafilters we obtain an order-preserving one-to-one correspondence  between finitely additive measures and ultrafilters on a given Boolean algebra.

\begin{corollary}
     For any Boolean algebra $\cB$,  $\langle \mathcal{UF}_{\!\cB}, \subseteq \rangle$ and $\langle \mathcal{FAM}_{0, 1, \cB}, \leq \rangle$ are order-isomorphic via the map $\Psi_{\cB} \colon \mathcal{UF}_{\cB} \to \mathcal{FAM}_{0, 1, \cB}$ defined by $\Psi_{\cB}(U) \coloneqq \Xi_{u}$ for every ultrafilter $U$ on $\cB$. 
\end{corollary}

In particular, ultrafilters are particular cases of finitely additive measures. 

\subsection{Notions of absoluteness}\label{n.2}

In this subsection, we introduce the basic elements of the notions of absoluteness. A reader unfamiliar with these notions may refer to~\cite{KunenFoM} and~\cite{Je2}.  We start by introducing the notion of \emph{relativization} of a formula to a class. 

\begin{definition}
    Let $\calM$ be a class and $\varphi$ a formula in the language of set theory. The \emph{relativization} of $\varphi$ to $\calM$, denoted by $\varphi^{\calM},$ is defined recursively based on the complexity of $\varphi$, as follows:

    \begin{enumerate}[label=(\arabic*)]
        \item If $\varphi$ is atomic, that is, of the form $x = y$ or $x \in y$, then $ \varphi^{\mathcal{M}} \coloneqq \varphi $.
        
        \item If $\varphi = \neg\psi,$ then $\varphi^{\mathcal{M}} \coloneqq \neg\psi^{\mathcal{M}}$.
        
        \item If $\varphi = \chi \wedge \psi, $ then $ \varphi^{\mathcal{M}} \coloneqq \chi^{\mathcal{M}} \wedge \psi^{\mathcal{M}}$. Similarly for the other logical connectives.
        
        \item If $\varphi = \exists x\psi(x), $ then $\varphi^{\mathcal{M}} \coloneqq \exists x\in \mathcal{M}[\psi^{\mathcal{M}}(x)]$. Similarly for the universal quantifier.
    \end{enumerate}
\end{definition}

We now introduce one of the central notions in this paper: \emph{absoluteness}. 

\begin{definition}
	Let $ \mathcal{M}$ and $\mathcal{N}$ be classes such that $\mathcal{M}\subseteq \mathcal{N}$ and $\varphi$ a formula in the language of set theory. Then: 
	\begin{enumerate}[label=(\arabic*)]
		\item We say that $\varphi$ is \emph{absolute for $ \mathcal{M}, \mathcal{N}$}, denoted by $ \mathcal{M}\preccurlyeq_{\varphi} \mathcal{N}, $ if $ \varphi(x_{1}, ..., x_{n}) $ is a formula with at most the free variables $ x_{1}, ..., x_{n}$, and 
        $$
        \forall \vec{a} \in \mathcal{M}^{n}[\varphi^{\mathcal{N}}(\vec{a}) \Leftrightarrow  \varphi^{\mathcal{M}}(\vec{a})],
        $$ where $\vec{a} = (a_{1}, \dots, a_{n})$ denotes an $n$-tuple of elements in $\mathcal{M}.$ 
		
		\item $\varphi$ is said to be \emph{absolute for $\mathcal{M}$} if it is absolute for $\mathcal{M}, V$, where $V$ is the universe of set theory.   
	\end{enumerate}
\end{definition}



Recall that a class $\calM$ is \emph{transitive} if, for any $x \in \calM$ and $y \in x$, it follows that $y \in \calM$. A quantifier of the form ``$\exists y \in x$'' or ``$\forall y \in x$'' is referred to as a \emph{bounded quantifier}, and a formula in which all quantifiers are bounded is called a \emph{$\Delta_{0}$ formula}. These formulas can be constructed recursively  and are well-known to be absolute for transitive classes (see~\cite{KunenFoM}). As a result, several elementary set-theoretic notions —such as the empty set, being a subset, transitive set, function,  bijective function,  finite set,  natural, rational, or real number, and basic operations like union, intersection, or Cartesian product— are absolute for transitive classes. On the other hand, we also have that the sets of natural and rational numbers are absolute for transitive models of $\ZFC$ as well, that is, if $\calM$ is a transitive class, then $\bbN^{\calM} \coloneqq \calM \cap \bbN = \bbN$ and $\bbQ^{\calM} \coloneqq \calM \cap \bbQ = \bbQ$. Somewhat more complex notions, such as being  an upper bound of a subset of real numbers, are also absolute. Using this, we can prove that the notions of \emph{supremum} and \emph{infimum} for subsets of real numbers are also absolute. Specifically,

\begin{lemma}\label{t96}
    Let $M, N$ be transitive models of $\ZFC$ such that $M\subseteq N.$ If $X \in M$ and  $X \subseteq \bbR^{M} \coloneqq \bbR \cap M,$ then $\sup^{M}(X) = \sup^{N}(X)$ and $\inf^{M}(X) = \inf^{N}(X)$.
\end{lemma}

\begin{PROOF}{\ref{t96}}
    Let $X \in M$ such that $X \subseteq \bbR^{M}.$ It is clear that, in $N,$ $\sup^{M}(X)$ is a upper bound of $X,$ so $\sup^{N}(X) \leq \sup^{M}(X).$ Now, towards a contradiction, working in $N,$ assume that $\sup^{N}(X) < \sup^{M}(X)$ and let $r \in \bbQ$ such that $\sup^{N}(X) < r < \sup^{M}(X).$ 

    Now, working in $M,$ since the notions of \emph{upper bound} and \emph{rational number} are absolute for transitive models, we have that $r$ is a rational number in $M$, and it is an upper bound of $X,$ hence $\sup^{M}(X) \leq r.$  This implies that $r < r$ in $N,$ a contradiction. Therefore $\sup^{N}(X) = \sup^{M}(X).$

    Finally, the proof for $\inf^{N}(X) = \inf^{M}(X)$ follows similar lines. 
\end{PROOF}

We focus now on the notion of real number. Although the sets of natural and rational numbers are absolute for transitive models, this is not the case for the set of real numbers. It is well known that this set is not absolute for transitive models of $\ZFC$. For instance, the forcing method 
can be used to construct, given a model $M$ of $\ZFC$, an extension $N$ of $M$ such that $\mathbb{R}^{M} \coloneqq \bbR \cap M \subsetneq \mathbb{R}^{N} \coloneqq N \cap \bbR$, that is, an extension of $M$ containing new real numbers not present in $M$. Similarly for the case of $\bbR^{n}$.  As a consequence, the following result holds in general. 

\begin{fact}\label{da5}
    If $M, N$ are transitive models of $\ZFC$ and $M \subseteq N$, then $(\bbR^{n})^{M} \subseteq (\bbR^{n})^{N}$.  In particular, if $a, b \in \bbR^{n} \cap M$ and $a \leq b$ then $[a, b]^{M} \coloneqq [a, b] \cap M \subseteq [a, b]^{N} \coloneqq [a, b] \cap N$ and $[a, b]^{M}$ is dense in $[a, b]^{N}$. 
\end{fact}

In many interesting cases inclusions in~\autoref{da5} are strict: $(\bbR^{n})^{M} \subsetneq (\bbR^{n})^{M}$ and similarly $[a, b]^{M} \subsetneq [a, b]^{N}$. 

\section{Riemann integration on Boolean algebras}\label{3}

The primary goal of this section is to prove the absoluteness of the Riemann integral (see~\autoref{ta} and~\autoref{s20}) and other related results concerning the absoluteness of integration over Boolean algebras (see~\autoref{3.3}). To achieve this, we will begin by extending the classical notion of  Riemann integral within the framework of Boolean algebras, as outlined in~\cite{CMUP} and~\cite[Ch.~3]{uribethesis}.

\subsection{Integration on Boolean algebras}\label{3.1.i}

In this subsection, fix a Boolean sub-algebra $\cB$ of $\calP(X)$ for some non-empty set $X$, and some finitely additive measure $\Xi \colon \cB \to [0, \delta],$ where $\delta$ is a positive real number.  

We start by defining partitions and their refinements: 

\begin{definition}\label{t5}
    \ 
    
    \begin{enumerate}[label=\normalfont(\arabic*)]
    
        \item $\Pbf^\Xi$ is the set of finite partitions of $X$ into sets in $\dom(\Xi)$.

        \item If $P, Q \in \bfP^{\Xi}$, we say that \emph{$Q$ is a refinement of $P$}, denoted by $Q \ll P$, if every element of $P$ can be finitely partitioned into elements of $Q.$ \index{$\ll$} \index{refinement}

        \item If $P = \langle P_{n} \colon n < n^{\ast} \rangle$ and $Q = \langle Q_{m} \colon m < m^{\ast} \rangle$ are in $\bfP^{\Xi},$ we define: $$P \sqcap Q \coloneqq \bigcup \{ P_{n} \cap Q_{m} \colon n < n^{\ast} \wedge m < m^{\ast} \} .$$ \index{$\sqcap$}
        \begin{figure}[H]
        \centering
        \begin{tikzpicture}[scale=0.7]
            \draw (0, 0) -- (0, 4);
            \draw (0, 0) -- (4, 0);
            \draw (4, 0) -- (4, 4);
            \draw (0, 4) -- (4, 4);

            \draw (6, 0) -- (6, 4);
            \draw (6, 0) -- (10, 0);
            \draw (10, 0) -- (10, 4);
            \draw (6, 4) -- (10, 4);

            \draw (12, 0) -- (12, 4);
            \draw (12, 0) -- (16, 0);
            \draw (16, 0) -- (16, 4);
            \draw (12, 4) -- (16, 4);     

            \draw[redun] (0, 1) -- (4, 1);
            \draw[redun] (0, 2) -- (4, 2);
            \draw[redun] (0, 3) -- (4, 3);
            \draw[redun] (1, 0) -- (1, 4);
            \draw[redun] (2, 0) -- (2, 4);
            \draw[redun] (3, 0) -- (3, 4);

            \draw[greenun] (6, 0.5) -- (10, 0.5);
            \draw[greenun] (6, 1.5) -- (10, 1.5);
            \draw[greenun] (6, 2.5) -- (10, 2.5);
            \draw[greenun] (6, 3.5) -- (10, 3.5);
            
            \draw[greenun] (6.5, 0) -- (6.5, 4);
            \draw[greenun] (7.5, 0) -- (7.5, 4);
            \draw[greenun] (8.5, 0) -- (8.5, 4);
            \draw[greenun] (9.5, 0) -- (9.5, 4);

            \draw[redun] (12, 1) -- (16, 1);
            \draw[redun] (12, 2) -- (16, 2);
            \draw[redun] (12, 3) -- (16, 3);
            \draw[redun] (13, 0) -- (13, 4);
            \draw[redun] (14, 0) -- (14, 4);
            \draw[redun] (15, 0) -- (15, 4);

            \draw[greenun] (12, 0.5) -- (16, 0.5);
            \draw[greenun] (12, 1.5) -- (16, 1.5);
            \draw[greenun] (12, 2.5) -- (16, 2.5);
            \draw[greenun] (12, 3.5) -- (16, 3.5);
            
            \draw[greenun] (12.5, 0) -- (12.5, 4);
            \draw[greenun] (13.5, 0) -- (13.5, 4);
            \draw[greenun] (14.5, 0) -- (14.5, 4);
            \draw[greenun] (15.5, 0) -- (15.5, 4);

            \node at (-0.5, 4) {$X$};
            \node at (5.5, 4) {$X$};
            \node at (11.5, 4) {$X$};

            \node[redun] at (4.4,0) {$P$};
            \node[greenun] at (10.4, 0) {$Q$};
            \node at (17, 0) {$P \sqcap Q$};
        \end{tikzpicture}
        
        \caption{A graphic example of $P \sqcap Q.$}
        \label{f46}
    \end{figure}
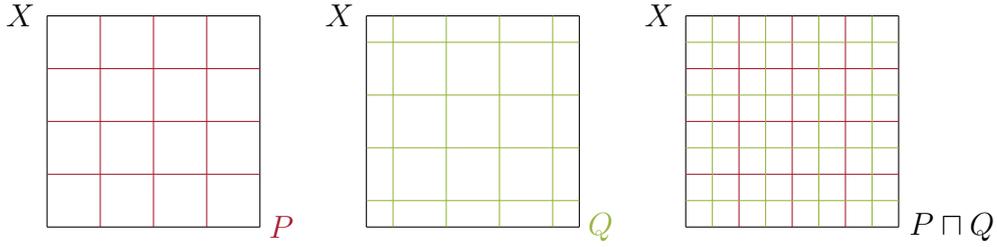    
    \end{enumerate}
\end{definition}

For example, it is clear that $\{ X \} \in \bfP^{\Xi}$ and if $P \in \bfP^{\Xi},$ then $P \ll \{ X \}$ and $P \ll P.$ Moreover, $\ll$ is a partial order on $\bfP^{\Xi}.$ Furthermore, for $P, Q \in \bfP^{\Xi}$, $P \sqcap Q$ is a common refinement of $P$ and $Q$:

\begin{lemma}\label{t9}
    If $P, Q \in \bfP^{\Xi},$ then $P \sqcap Q \in \bfP^{\Xi}$ and $P \sqcap Q \ll P, Q.$
\end{lemma}


We can now define the integral with respect to $\Xi$:

\begin{definition}\label{t14}
  \index{$\supsum(f, P)$} \index{$\infsum(f, P)$}
  
  Let $f\colon X\to\bbR$ be a bounded function. We define:
  
  \begin{enumerate}   [label=\normalfont(\arabic*)]   
      \item For any $P\in\Pbf^\Xi$, 
      \[\supsum_{\Xi}(f,P) \coloneqq \sum_{b\in P}\sup(f[b])\Xi(b)\text{\ and \ }\infsum_{\Xi}(f,P) \coloneqq \sum_{b\in P}\inf(f[b])\Xi(b).\] 
      
      \item $\overline{\int_{X}} fd\Xi \coloneqq \inf \left \{\supsum(f,P) \colon  P\in\Pbf^\Xi \right \}$  and 
      $\underline{\int_{X}}fd\Xi \coloneqq \sup \left \{\infsum(f,P) \colon  P\in\Pbf^\Xi \right \}$. \index{$\overline{\int}_{X} d \Xi$} \index{$\underline{\int}_{X} f d \Xi$}
      
      \item We say that $f$ is \emph{$\Xi$-integrable on $X$}, denoted by $f \in \cI(\Xi)$, iff $\overline{\int_{X}}fd\Xi=\underline{\int_{X}}fd\Xi$. In this case, this common value is denoted  by $\int_{X} fd\Xi$. 
  \end{enumerate}

    Naturally, when the context is clear, we omit the superscript ``$\Xi$'' in ``$\supsum^{\Xi}(f, P)$'' and ``$\infsum^{\Xi}(f, P)$''. 
\end{definition}

For example, it is clear that any constant function is $\Xi$-integrable. Concretely, if for all $x \in X, $ $ f(x) = c \in \bbR,$ then $\int_{X} f(x) \, d \Xi = c \Xi(X)$.

\begin{lemma}\label{t18}
    Let $f \colon X \to \bbR$ be a bounded function. If $P, Q \in \bfP^{\Xi}$ and $Q \ll P,$ then $$
    \infsum(f, P) \leq \infsum(f, Q) \leq \supsum(f, Q) \leq \supsum(f, P).
    $$ 
    As a consequence, $\supsum(f, Q) - \infsum(f, Q) \leq \supsum(f, P) - \infsum(f, P).$
\end{lemma}

We can use $P \sqcap Q$ to get the following result. 

\begin{corollary}\label{t22}
    If $P, Q \in \bfP^{\Xi},$ then $\infsum(f, P) \leq \supsum(f, Q).$ As a consequence, the following inequality holds:  $$
    \underline{\int_{X}} f d  \Xi \leq \overline{\int_{X}} f d  \Xi.
    $$
\end{corollary}


We will now proceed to prove the result we refer to as the \emph{Criterion of $\Xi$-Integrability}, which, analogous to the case of the Riemann integral, allows us to characterize the $\Xi$-integrability of a function in terms of the existence of suitable partitions. 

\begin{theorem}\label{t29}
    Let $f \colon X \to \bbR$ be a bounded function. Then, $f$ is $\Xi$-integrable if, and only if, for all $\varp > 0,$ there exists a partition $P \in \bfP^{\Xi}$ such that $ \supsum(f, P) - \infsum(f, P) < \varp.$
\end{theorem}

\begin{PROOF}{\ref{t29}}
    On the one hand, assume that $f \in \cI(\Xi)$ and let $\varp > 0.$ By basic properties of $\sup$ and $\inf,$ we can find $P, Q \in \bfP^{\Xi}$ such that $$ \int_{X} f d \Xi - \frac{\varp}{2} < \infsum(f, P) \text{ and } \supsum(f, Q) < \int_{X} f d \Xi + \frac{\varp}{2}. $$ Consider $R \coloneqq P \sqcap Q.$ By~\autoref{t9}, $R \in \bfP^{\Xi}$ and it is a common refinement of $P$ and $Q$. So, by virtue of  \autoref{t18}, $\infsum(f, P) \leq \infsum(f, R) \text{ and } \supsum(f, Q) \leq \supsum(f, R).$ Therefore, 
    $$ \int_{X} f d \Xi - \frac{\varp}{2} < \infsum(f, R) \text{ and } \supsum(f, R) < \int_{X} f d \Xi + \frac{\varp}{2}. 
    $$ 
    Thus, $\supsum(f, R) - \infsum(f, R) < \varp.$ 

    On the other hand, let $P \in \bfP^{\Xi}$ such that $\supsum(f, P) - \infsum(f, P) < \varp.$ Hence, 
    $$
    \overline{\int_{X}} f d \Xi \leq \supsum(f, P) < \infsum(f, P) + \varp \leq \underline{\int_{X}} f d \Xi + \varp.
    $$ Since $\varp$ is arbitrary, by~\autoref{t22} it follows that $f \in \cI(\Xi).$ 
\end{PROOF}

In fundamental aspects, the integral with respect to finitely additive measures behaves similarly to the Riemann integral, that is, we have available the basic properties of the integral such as those presented in \cite[Sec.~3.5]{uribethesis} and~\cite{CMUP}, for instance: 

\begin{lemma}
    Let $f, g \colon X \to \bbR$ be $\Xi$-integrable functions and $c \in \bbR.$ Then: 

    \begin{enumerate}[label = \normalfont (\alph*)]
        \item\label{t61} $fg$ and $cf$ are $\Xi$-integrable and  $\int_{X} (cf) d \Xi = c \int_{X} f d \Xi$.

        \item \label{t50} Let $\{ f_{i} \colon i < n \}$ be a finite sequence of $\Xi$-integrable functions. Then $\sum_{i < n} f_{i}$ is $\Xi$-integrable and $$\int_{X} \left( \sum_{i < n} f_{i}\right) d \Xi = \sum_{i < n} \left( \int_{X} f_{i}  d \Xi \right) \! .$$

        \item If $f \leq g$ then $ \int_{X} f d \Xi \leq \int_{X} g d \Xi.$
    \end{enumerate}
\end{lemma}



    


\subsection{Riemann integration on $\bbR^{n}$}\label{3.2.r}

In this Subsection, on the one hand, following the ideas of~\cite{Munkres}, we review some basic notions of Riemann integration over rectangles in $\bbR^{n}$. On the other hand, we show that this is a particular case of the integration introduced in~\autoref{3.1.i}. 

For the rest of this article, we fix $n \in \bbN$ and two points $a = (a_{0}, \dots, a_{n-1})$, $b = (b_{0}, \dots, b_{n-1})$ in $\bbR^{n}$ such that $a \leq b$, that is, for any $i < n$, $a_{i} \leq b_{i}$. We also denote by $\Lambda^{n}$ the Lebesgue measure on $\bbR^{n}$ and $[a, b] \coloneqq \prod_{i < n}[a_{i}, b_{i}]$. A \emph{partition} of some interval $[x, y] \subseteq \bbR$ is a finite collection $P = \{ p_{i} \colon i < m \}$ of increasing real numbers in $[x, y]$ such that $p_{0} = x$ and $p_{m-1} = y$. Each interval $[p_{i}, p_{i+1}]$ is called a \emph{sub-interval} determined by $P$. Similarly, a \emph{partition} of the rectangle $[a, b]$ is a $n$-tuple $\langle P_{0}, \dots, P_{n-1} \rangle$ such that for any $i < n$, $P_{i}$ is a partition of the interval $[a_{i}, b_{i}]$. If for any $i < n$, $I_{i}$ is some sub-interval determined by $P_{i}$, then $R \coloneqq \prod_{i < n} I_{i}$ is a \emph{sub-rectangle of $[a, b]$ determined by $P$}. Denote by $\calS_{P}$ the collection of all sub-rectangles of $[a, b]$ determined by $P$. 

If $f \colon [a, b] \to \bbR$ is a bounded function, for $R \in \calS_{P}$, define $m_{R}(f) \coloneqq \inf \{ f(x) \colon x \in R \}$ and $M_{R}(f) \coloneqq \sup \{ f(x) \colon x \in R \}$. We also define 
$$
L(f, P) \coloneqq \sum_{R \in \calS_{P}} m_{R}(f) \Lambda^{n}(R) \text{\ and \ } U(f, P) \coloneqq \sum_{R \in \calS_{P}} M_{R}(f) \Lambda^{n}(R).
$$

We say that $f$ is \emph{Riemann integrable} iff for any $\varp > 0$, there exists some partition $P$ of $[a, b]$ such that 
$ U(f, P) - L(f, P) < \varp.$
In this case, we define the Riemann integral of $f$ on $[a, b]$, as follows:
\begin{equation*}
    \begin{split}
        \int_{[a, b]} fd \Lambda^{n} & \coloneqq \sup \left\{ L(f, P) \colon P \text{ is a partition of } [a, b] \right\}.
    \end{split}
\end{equation*}

Recall that, a function $\sigma \colon [a, b] \to \bbR$ is a \emph{step function} if there exists some partition $P = \langle P_{j} \colon j < m \rangle$ of $[a, b]$ such that $\sigma$ is constant on the interior of each $R \in \calS_{P}$. Notice that the value of $\sigma$ on the boundary of each $R$ is not important for the purposes of integration, since this boundary has measure zero in $\mathbb{R}^n$. In the following result, we can characterize Riemann integration in terms of step functions (see~\cite[Sec.~$7$.A.$3$)]{Jones}).

\begin{lemma}
    A bounded function $f \colon [a, b] \to \bbR$ is Riemann integrable iff for any $\varp > 0$ there are step functions $\sigma, \tau$ on $[a, b]$ such that $\sigma \leq f \leq \tau$, except on some measure zero set, and 
    $$
    \int_{[a, b]} (\tau - \sigma)d \Lambda^{n} < \varp.
    $$ 
\end{lemma}

To show that the Riemann integral is a particular case of the integration defined previously in~\autoref{3.1.i}, we introduce the following notation.

\begin{definition}\label{t24}
    Let $E\subseteq X$. We define  $\cB|_E \coloneqq \{ E \cap b \colon b \in \cB \}$, which is a Boolean sub-algebra of $\calP(E)$. When $E\in\cB$, we denote $\Xi|_E \coloneqq \Xi \rest \cB|_E$,  which is a finitely additive measure on $\cB |_{E}.$
\end{definition} 

\begin{definition}
    Let $\cR^{n}$ be the collection of all subsets of $\bbR^{n}$ which are finite unions of rectangles of the form $[c,d) = \prod_{i < n}[c_{i}, d_{i})$  with $c = (c_{0}, \dots, c_{n-1})$, $d = (d_{0}, \dots, d_{n-1})$ and $c \leq d$ in $\bbR^{n}$.
\end{definition}

Notice that $\cR^{n}$ is not a Boolean algebra because, although $\emptyset \in \cR^{n}$ and it is closed under $\cup$, $\cap$, and set difference,
$\bbR^{n} \notin \cR$. 

\begin{fact}\label{r15}
    Define $\lambda^{n} \coloneqq \Lambda^{n} \rest \cR^{n}$. Then: 

    \begin{enumerate}[label = \normalfont (\alph*)]
        \item $\lambda^{n}(\emptyset) = 0$. 

        \item $\lambda^{n} \left( \bigcup_{m \in \bbN} I_{m} \right) = \sum_{m \in \bbN} \lambda^{n}(I_{m})$ whenever $\langle I_{m} \colon m \in \bbN \rangle$ is a pairwise disjoint sequence of elements in $\cR^{n}$. 

        \item For any $c, d \in \bbR^{n}$ with $c \leq d$, $\lambda^{n}([c, d)) = \prod_{i < n} (d_{i} - c_{i})$.   
    \end{enumerate}
\end{fact}


By~\autoref{t24}, $\cR^{n} |_{[a, b]}$ is a Boolean sub-algebra of $\calP([a, b])$ and $\lambda^{n} |_{[a, b]}$ is a finitely additive measure on $\calP([a, b])$. Using this, it follows that the Riemann integral over $\bbR^{n}$ is a particular case of our integral over Boolean algebras.

\begin{theorem}\label{r20}
    Let $f\colon [a,b] \to \bbR$ be a bounded function. Then, $f$ is Riemann integrable if, and only if, it is $\lambda^{n}|_{[a,b]}$-integrable. In this case, 
    $$
    \int_{[a, b]} f d \Lambda^{n} = \int_{[a, b]} f d \lambda^{n} |_{[a, b]}.
    $$
\end{theorem}

\begin{PROOF}{\ref{r20}}
    Straightforward using~\autoref{t29}. 
\end{PROOF}

\subsection{Absoluteness of the Riemann integral}\label{3.3}

In this subsection, we present several results related to the absoluteness of the integral that we defined in~\autoref{3.1.i}. In particular, we use~\autoref{r20} to prove the main result in this article: the Riemann integral is absolute for transitive models of $\ZFC$ (see~\autoref{ta} and~\autoref{s20}). 

We start by proving some preliminaries results in $\ZFC$. 

\begin{theorem}\label{t93}
    Let $\cB_{0}, \cB_{1} \subseteq \calP(X)$ be Boolean algebras such that $\cB_{0} \subseteq \cB_{1},$ and $\Xi_{0}, \Xi_{1}$ be finitely additive measures on $\cB_{0}$ and $ \cB_{1}$, respectively, such that $\Xi_{0} \subseteq \Xi_{1}.$ Let $f \colon X \to \bbR$ be a bounded function. Then $ f \in \cI(\Xi_{0}) $ implies $ f \in \cI(\Xi_{1})$. In this case, $$ \int_{X} f  d \Xi_{0} = \int_{X} f d \Xi_{1}.$$ The converse implication holds whenever $\cB_{0} = \cB_{1}.$ 
\end{theorem}

\begin{PROOF}{\ref{t93}}
    For any $P \in \bfP^{\Xi_{0}}$ we have that $\supsum_{\Xi_{0}}(f, P) = \supsum_{\Xi_{1}}(f, P)$ and $\infsum_{\Xi_{0}}(f, P)  =  \infsum_{\Xi_{1}}(f, P)$. Consequently,  
    $\{\supsum_{\Xi_{0}}(f, P) \colon P \in \bfP^{\Xi_{0}}\} \subseteq \{ \supsum_{\Xi_{1}}(f, P) \colon P \in \bfP^{\Xi_{1}} \}$ and, in a similar way,
    $ \{ \infsum_{\Xi_{0}}(f, P) \colon P \in \bfP^{\Xi_{0}} \} \subseteq \{ \infsum_{\Xi_{1}}(f, P) \colon P \in \bfP^{\Xi_{1}} \}.$ Therefore, $$
    \overline{\int_{X}} f d  \Xi_{1} = \inf \{ \supsum_{\Xi_{1}}(f, P) \colon P \in \bfP_{\Xi_{1}} \} \leq \inf \{ \supsum_{\Xi_{0}}(f, P) \colon P \in \bfP^{\Xi_{0}} \} = \overline{\int_{X}} f  d \Xi_{0}, \text{ and}
    $$ $$ \underline{\int_{X}} f d \Xi_{0}  = \sup \{ \infsum_{\Xi_{0}}(f, P) \colon P \in \bfP_{\Xi_{0}} \} \leq \sup \{ \infsum_{\Xi_{1}}(f, P) \colon P \in \bfP^{\Xi_{1}} \} = \underline{\int_{X}} f d  \Xi_{1}.
    $$ 
    As a consequence,
    \begin{equation}\label{e45}
        \text{$ \underline{\int_{X}} f d  \Xi_{0}  \leq \underline{\int_{X}} f  d \Xi_{1} \leq \overline{\int_{X}} f  d \Xi_{1} \leq \overline{\int_{X}} f d  \Xi_{0}.$}
    \end{equation}

    Thus, if $f$ is $\Xi_{0}$-integrable, then $f$ is $\Xi_{1}$-integrable, and it is clear that their values coincide.

    To prove the converse, notice that if $\cB_{0} = \cB_{1},$ then $\Xi_{0} = \Xi_{1}$ and $\bfP^{\Xi_{0}} = \bfP^{\Xi_{1}}$ and therefore, inequalities in~(\ref{e45}) ares really equalities. Thus, $f \in \cI(\Xi_{0})$ iff $f \in \cI(\Xi_{1}),$ and it is clear that the values of the integrals are the same. 
\end{PROOF}

We can now prove that the integral with respect to finitely additive measures in Boolean algebras is absolute for transitive models of $\ZFC$.

\begin{theorem}\label{t100}
    Let $M, N$ be transitive models of $\ZFC$ such that $M\subseteq N.$ Let $\cB, \Xi, X$ in $M$ be such that $\cB \subseteq \calP(X)$ is a Boolean algebra and $\Xi$ is a finitely additive measure on $\cB$. If $f \colon X \to \bbR$ is a bounded function in $M$, then 
    $$
    [f \in \cI(\Xi)]^{M} \Leftrightarrow [f \in \cI(\Xi)]^{N}.
    $$
    In this case, $$ \left(  \int_{X} f d\Xi  \right)^{  N} = \left( \int_{X} f d \Xi \right)^{M} \! \! \!. $$   
\end{theorem}

\begin{PROOF}{\ref{t100}}
    By \autoref{t96}, we have that 
    $$
    (\bfP^{\Xi})^{M} = (\bfP^{\Xi})^{N}, \ \supsum^{M}_{\Xi}(f, P) = \supsum_{\Xi}^{N}(f, P), \text{ and } \infsum^{M}_{\Xi}(f, P) = \infsum_{\Xi}^{N}(f, P).$$ Consequently,  
    $
    \{\supsum^{M}_{\Xi}(f, P) \colon P \in (\bfP^{\Xi})^{M}\} = \{ \supsum^{N}_{\Xi}(f, P) \colon P \in (\bfP^{\Xi})^{N} \},$ and in a similar way,
    $ \left\{\infsum^{M}_{\Xi}(f, P) \colon P \in (\bfP^{\Xi})^{M} \} = \{ \infsum^{M}_{\Xi}(f, P) \colon P \in (\bfP^{\Xi})^{N} \right\}. 
    $ Therefore,  by taking infimum in the first equality and supremum in the second, we get: 
    $$ 
    \left( \underline{\int_{X}} f d  \Xi \right)^{  M} = \left( \underline{\int_{X}} f d  \Xi \right)^{  N} \text{ and } \left( \overline{\int_{X}} f d  \Xi \right)^{  M}  = \left( \overline{\int_{X}} f d  \Xi \right)^{N} \! \! \!, 
    $$ which proves the result.
\end{PROOF}

\begin{corollary}\label{t102}
    Let $M, N$ be transitive models of $\ZFC$ such that $M \subseteq N$. Let $\Xi_{0}, \cB_{0} \in M $ and $\Xi_{1}, \cB_{1} \in N.$ Assume that $\Xi_{0}, \Xi_{1}$ are finitely additive measures on $\cB_{0} $ and $\cB_{1}$ respectively, such that $\Xi_{0} \subseteq \Xi_{1}$ and $\cB_{0} \subseteq \cB_{1} \subseteq \calP(X)$ for some set $X \in M.$ Let $f \colon X \to \bbR$ be a bounded function in $M$. Then: 
    $$
    [f \in \cI(\Xi_{0})]^{M} \Rightarrow [f \in \cI(\Xi_{1})]^{N}.
    $$ In this case, 
    $$ \left(  \int_{X} f d\Xi_{0}  \right)^{  M} = \left( \int_{X} f d \Xi_{1} \right)^{N} \! \! \!. 
    $$   
\end{corollary}

We now focus specifically on the Riemann integral on $\bbR^{n}$. Notice that the property of absoluteness for the Riemann integral does not follow directly from~\autoref{t100}: if $f$ is defined on $[a, b]$ in $M$, as established in~\autoref{n.2}, it is possible that $[a, b]^{M} \subsetneq [a, b]^{N}$, which implies that when we consider  $f$ in $N$, it could be not defined on the whole $[a, b]^{N}$, according to the interpretation of $[a, b]^{N}$. Therefore, in $N$, it is necessary to extend the function and integrate over a larger set. Due to this situation, the proof for the Riemann integral requires additional effort.

Recall from~\autoref{s16} and~\autoref{ss0} that, if $h \colon X \to Y$ is a function and $\cB$ is a Boolean sub-algebra of $\calP(X)$, then the collection $h^{\to}(\cB) \coloneqq \{ A \subseteq Y \colon h^{-1}[A] \in \cB \}$ is a Boolean sub-algebra of $\calP(Y)$ and, if $\Xi$ is a finitely additive measure on $\cB$, then $\Xi_{h}$ is a finitely additive measure on $h^{\to}(\cB)$. 

\begin{lemma}\label{s5}
    Let $X$, $Y$ be non-empty sets, $h \colon X \to Y$ a function, $\cB_{X}$ a Boolean sub-algebra of $\calP(X)$, $\cB_{Y} \coloneqq h^{\to}(\cB_{X})$, and $\Xi$ a finitely additive measure defined on $\cB_{X}$. Let $f \colon Y \to \bbR$ be a bounded function. Then:

    \begin{enumerate}[label = \normalfont (\alph*)]
        \item\label{s5.a} For any $P \in \bfP^{\Xi_{h}}$, there are $P^{\bullet} \in \bfP^{\Xi_{h}}$ and $Q \in \bfP^{\Xi}$ such that $P^{\bullet} \ll P$ and $$ \supsum_{\Xi_{h}}(f, P^{\bullet}) = \supsum_{\Xi}(f \circ h, Q) \text{ and }  \infsum_{\Xi_{h}}(f, P^{\bullet}) = \infsum_{\Xi}(f \circ h, Q).$$

        \item\label{s5.b} If $h$ is one-to-one, then for any $Q \in \bfP^{\Xi}$ there exists some $P \in \bfP^{\Xi_{h}}$ such that $$\supsum_{\Xi}(f \circ h, Q) = \supsum_{\Xi_{h}}(f, P) \text{ and } \infsum_{\Xi}(f \circ h, Q) = \infsum_{\Xi_{h}}(f, P).$$ 
    \end{enumerate}
\end{lemma}

\begin{PROOF}{\ref{s5}}
    \ref{s5.a}: Let $P \in \bfP^{\Xi_{h}}$. Define $P^{\bullet} \coloneqq P \sqcap \{ \ran h, (\ran h)^{c} \}$. Hence, $P^{\bullet} \in \bfP^{\Xi_{h}}$, it is a refinement of $P$, and for any $A \in P^{\bullet}$, either $A \subseteq \ran h$ or $A \cap \ran h = \emptyset$. Consider the set $R \coloneqq \{ A \in P^{\bullet} \colon A \subseteq \ran h \}$. Notice that, if $A \in R$ then $A = h[h^{-1}[A]]$ and, otherwise, $\Xi_{h}(A) = \emptyset$. Define $Q \coloneqq \{ h^{-1}[A] \colon A \in R \}$. Clearly $Q \in \bfP^{\Xi}$. As a consequence,  
    \begin{equation*}
        \begin{split}
            \supsum_{\Xi_{h}}(f, P^{\bullet}) & = \sum_{A \in P^{\bullet}} \sup(f[A]) \Xi_{h}(A) = \sum_{A \in R} \sup(f[A]) \Xi_{h}(A) + \sum_{A \in P^{\bullet} \setminus R} \sup(f[A]) \Xi_{h}(A)\\
            & = \sum_{A \in R} \sup(f \left[h[h^{-1}[A]] \right]) \Xi(h^{-1}[A]) = \sum_{B \in Q} \sup(f \circ h[B]) \Xi(B)\\
            & = \supsum_{\Xi}(f \circ h, Q). 
        \end{split}
    \end{equation*}

    Similarly, $ \infsum_{\Xi_{h}}(f, P^{\bullet}) = \infsum_{\Xi}(f \circ h, Q)$. 

    \ref{s5.b}: Assume that $h$ is one-to-tone. Let $Q \in \bfP^{\Xi}$.  Define $P \coloneqq \{ h[B] \colon B \in Q \} \cup \{ (\ran h)^{c} \}$. Clearly, $P \in \bfP^{\Xi_{h}}$, and
    \begin{equation*}
        \begin{split}
            \supsum_{\Xi_{h}}(f, P) & = \sum_{A \in P} \sup(f[A]) \Xi_{h}(A)  = \sum_{B \in Q} \sup(f[h[B]]) \Xi_{h}(h[B])\\
            & = \sum_{B \in Q} \sup(f \circ h [B]) \Xi(B) = \supsum_{\Xi}(f \circ h, Q). 
        \end{split}
    \end{equation*}
    The proof for the lower sum is similar. 
\end{PROOF}

As a consequence, under the conditions  of~\autoref{s5}, integrability is preserved under composition as well as the value of the integrals. 

\begin{theorem}\label{s2}
    Let $X$, $Y$ be non-empty sets, $h \colon X \to Y$ a function, $\cB_{X}$ a Boolean sub-algebra of $\calP(X)$, $\cB_{Y} \coloneqq h^{\to}(\cB_{X})$, and $\Xi$ a finitely additive measure on $\cB_{X}$. Let $f \colon Y \to \bbR$ be a bounded function. Then: 
    
    \begin{enumerate}[label = \normalfont (\alph*)]
        \item\label{s2.1} If $f \in \cI(\Xi_{h})$ then $f \circ h \in \cI(\Xi)$ and 
         $$
        \int_{Y} f d \Xi_{h} = \int_{X} f \circ h d \Xi. 
        $$

        \item\label{s2.2} If $h$ is one-to-one, then $f \circ h \in \cI(\Xi)$ implies $f \in \cI(\Xi_{h})$ and 
         $$  \int_{X} f \circ h d \Xi = 
        \int_{Y} f d \Xi_{h}. 
        $$
    \end{enumerate}
\end{theorem}

\begin{PROOF}{\ref{s2}}
    \ref{s2.1}: Assume that $f \in \cI(\Xi_{h})$ and $\varp > 0$. So we can find a partition $P \in \bfP^{\Xi_{h}}$ such that $\supsum_{\Xi_{h}}(f, P) - \infsum_{\Xi_{h}}(f, P) < \varp$. Let $P^{\bullet} \in \bfP^{\Xi_{h}}$ and $Q \in \bfP^{\Xi}$ as in~\autoref{s5}~\ref{s5.a}. Since $P^{\bullet}$ is a refinement of $P$, by~\autoref{t18}, 
    $$
    \supsum_{\Xi}(f \circ h, Q) -  \infsum_{\Xi}(f \circ h, Q) = \supsum_{\Xi_{h}}(f, P^{\bullet}) - \infsum_{\Xi_{h}}(f, P^{\bullet}) \leq \supsum_{\Xi_{h}}(f, P) - \infsum_{\Xi_{h}}(f, P) < \varp.
    $$ 
    This shows that $f \circ h$ is $\Xi_{h}$-integrable by applying~\autoref{t29}. 

     We now deal with the value of the integral. On the one hand, let $P \in \bfP^{\Xi_{h}}$. By applying~\autoref{s5}~\ref{s5.a}, we can find $P^{\bullet} \ll P$  in $\bfP^{\Xi_{h}}$ and $Q \in \bfP^{\Xi}$ such that 
    $$
       \infsum_{\Xi_{h}}(f, P) = \infsum_{\Xi}(f \circ h, Q) \leq  \int_{X} f \circ h d \Xi \leq  \supsum_{\Xi}(f \circ h, Q) = \supsum_{\Xi_{h}}(f, P^{\bullet}) \leq \supsum_{\Xi_{h}}(f, P).
    $$

    As a consequence, $\int_{X} f \circ h d \Xi = \int_{Y} f d \Xi_{h}$.
    
    \ref{s2.2}: Assume that $h$ is one-to-one, $f \circ h$ is $\Xi$-integrable, and let $\varp > 0$. By~\autoref{t29}, there exists some $Q \in \bfP^{\Xi}$ such that $\supsum_{\Xi}(f \circ h, Q) - \infsum_{\Xi}(f \circ h, Q) < \varp$. Consider $P \in \bfP^{\Xi_{h}}$ as in~\autoref{s5}~\ref{s5.b}. As a consequence, 
    $$
        \supsum_{\Xi_{h}}(f, P) - \infsum_{\Xi_{h}}(f, P) = \supsum_{\Xi}(f \circ h, Q) - \infsum_{\Xi}(f \circ h, Q) < \varp.
    $$
    Thus, by~\autoref{t29}, $f \in \cI(\Xi_{h})$. Finally, the value of the integral follows by applying~\ref{s2.1}. 
\end{PROOF}

\begin{corollary}\label{da1}
    Let $X, Y$ be non-empty sets such that $X \subseteq Y$, $\cB_{X}$, $\cB_{Y}$ Boolean algebras on $\calP(X)$ and $\calP(Y)$, respectively, $\Xi^{X}$, $\Xi^{Y}$ finitely additive measures on $\cB_{X}$  and $\cB_{Y}$, respectively, and $g \colon Y \to \bbR$ a bounded function. Assume that $\cB_{Y} \subseteq h^{\to}(\cB_{X})$ and $\Xi_{Y} \subseteq \Xi_{h}^{X}$. If $g \in \cI(\Xi^{Y})$ then $g \rest X \in \cI(\Xi^{X})$ and 
    \begin{equation}\label{e66}
        \text{$\int_{Y} g d \Xi^{Y} = \int_{X} g \rest X d \Xi^{X}.$}
    \end{equation}
\end{corollary}

\begin{PROOF}{\ref{da1}}
    Let $h \colon X \to Y$ be the inclusion function and assume that $g \in \cI(\Xi^{Y})$. Since $\cB_{Y} \subseteq h^{\to}(\cB_{X})$ and $\Xi^{Y} \subseteq \Xi^{X}_{h}$, by applying~\autoref{t93} it follows that $g \in \cI(\Xi_{h}^{X})$ and 
    \begin{equation}\label{e67}
        \text{$\int_{Y} g d \Xi^{Y} = \int_{Y} g d \Xi_{h}^{X}.$}
    \end{equation}
    On the other hand, by~\autoref{s2}, we have that $g \rest X = g \circ h \in \cI(\Xi^{X})$ and
    \begin{equation}\label{e60}
        \text{$\int_{Y} g d \Xi_{h}^{Y} = \int_{X} g \rest X d \Xi^{X}.$}
    \end{equation}
    Finally,~(\ref{e66}) follows from~(\ref{e67}) and~(\ref{e60}). 
\end{PROOF}

We are finally ready to prove~\autoref{ta}, the main result of this paper. 

\begin{theorem}\label{s20}
    Let $M$, $N$ be transitive models of $\ZFC$ such that $M \subseteq N$, $n \in \bbN$, and  $a, b \in \bbR^{n} \cap M$ with $a \leq b$. In $M$, let $f$ be a real-valued function on $[a, b]$. Then, the following statements are equivalent:  

    \begin{enumerate}[label = \normalfont (\roman*)]
        \item\label{s20.1} $f$ is Riemann integrable in $M$.  

        \item\label{s20.2} In $N$, there exists some Riemann integrable function $g \colon [a, b] \to \bbR$ extending $f$. 
    \end{enumerate}

    If either~\ref{s20.1} or~\ref{s20.2} holds, then: 
    \begin{equation}\label{e-s20.1}
        \text{$
    \left( \int_{[a, b]} f d \lambda^{n} |_{[a, b]} \right)^{M} = \left(  \int_{[a, b]} g d \lambda^{n} |_{[a, b]} \right)^{N}\! \! \!.
    $}
    \end{equation}

    Moreover, the function $g$  in~\ref{s20.2} is unique except in a Lebesgue measure zero set: if $g^{\ast}$ is another Riemann integrable function on $[a, b]$ extending $f$, then there exists some Lebesgue measure zero set $E \subseteq [a, b]$ such that, for any $x \in [a, b] \setminus E$, $g^{\ast}(x) = g(x)$. 
\end{theorem}

 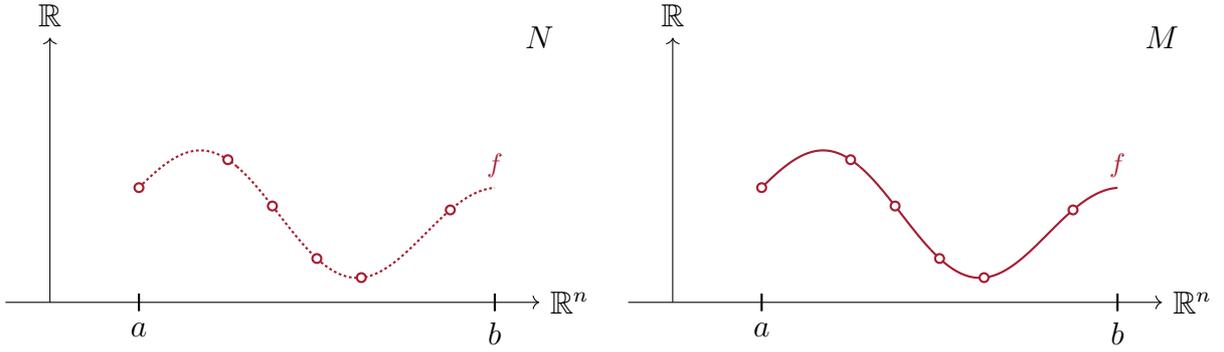
\begin{figure}
        \centering
        \begin{tikzpicture}[scale=1.17] 
    
        \draw[->] (-0.5, 0) -- (5.5, 0) node[right] {$\bbR^{n}$};
        
        \draw[->] (0, 0) -- (0, 3) node[above] {$\bbR$};
    
        \node at (5.5, 3) {$N$};
    
        \draw[thick] (1, 0.1) -- (1, -0.1) node[below] {$a$};
        \draw[thick] (5, 0.1) -- (5, -0.1) node[below] {$b$};
       \node[redun] at (5, 1.55) {\footnotesize $f$};

    
        \draw[domain=1:2.1, smooth, thick, redun, dash pattern={on 1pt off 1pt on 1pt off 1pt}]  plot (\x, {1 + 0.5*sin(2*(\x-1) r) + 0.3*cos((\x-1) r)});
            
        \draw[domain=2.1:4, smooth, thick, redun, dash pattern={on 1pt off 1pt on 1pt off 1pt}]  plot (\x, {1 + 0.5*sin(2*(\x-1) r) + 0.3*cos((\x-1) r)});
            
        \draw[domain=4:5, smooth, thick, redun, dash pattern={on 1pt off 1pt on 1pt off 1pt}] plot (\x, {1 + 0.5*sin(2*(\x-1) r) + 0.3*cos((\x-1) r)});
    
        \foreach \x in {1, 2, 2.5, 3,  3.5, 4.5} {
            \fill[white] (\x, {1 + 0.5*sin(2*(\x-1) r) + 0.3*cos((\x-1) r)}) circle(0.05);
            
            \draw[thick, redun] (\x, {1 + 0.5*sin(2*(\x-1) r) + 0.3*cos((\x-1) r)}) circle(0.05);
        }

        \draw[->] (6.5, 0) -- (12.5, 0) node[right] {$\bbR^{n}$};
        
        \draw[->] (7, 0) -- (7, 3) node[above] {$\bbR$};
    
        \node at (12.5, 3) {$M$};
    
        \draw[thick] (8, 0.1) -- (8, -0.1) node[below] {$a$};
        \draw[thick] (12, 0.1) -- (12, -0.1) node[below] {$b$};
    
    
        \draw[domain=8:9.1, smooth, thick, redun]  plot (\x, {1 + 0.5*sin(2*(\x-8) r) + 0.3*cos((\x-8) r)});
        
        \draw[domain=9.1:9.9, smooth, thick, redun]  plot (\x, {1 + 0.5*sin(2*(\x-8) r) + 0.3*cos((\x-8) r)});
        
        \draw[domain=9.9:10.5, smooth, thick, redun]  plot (\x, {1 + 0.5*sin(2*(\x-8) r) + 0.3*cos((\x-8) r)});
        
        \draw[domain=10.5:12, smooth, thick, redun]  plot (\x, {1 + 0.5*sin(2*(\x-8) r) + 0.3*cos((\x-8) r)});
    
        \foreach \x in {8, 9, 9.5, 10, 10.5, 11.5} {
            \fill[white] (\x, {1 + 0.5*sin(2*(\x-8) r) + 0.3*cos((\x-8) r)}) circle(0.05);
            \draw[thick, redun] (\x, {1 + 0.5*sin(2*(\x-8) r) + 0.3*cos((\x-8) r)}) circle(0.05);
        }
    
           \node[redun] at (12, 1.55) {\footnotesize $f$};
    \end{tikzpicture}
    
        \caption{The graph on the right represents the situation from the perspective of $M$: if $f$ is Riemann integrable, then it is a continuous function except of a set of measure zero. On the other hand, the graph on the left illustrates the situation from the perspective of $N$: with the appearance of new real numbers in $[a, b]$, the function $f$ is no longer defined over the entire rectangle $[a, b]$.}
        \label{fig:enter-label}
    \end{figure}

\begin{PROOF}{\ref{s20}}
    \ref{s20.2} $\Rightarrow$~\ref{s20.1}:  Working in $N$, assume that $g \colon [a, b] \to \bbR$ is a Riemann integrable function extending $f$.  With the intention of applying~\autoref{da1}, define $X \coloneqq [a, b]^{M}$, $Y \coloneqq [a, b]$, $h \colon X \to Y$ the inclusion map, $\cB_{X} \coloneqq \cR^{n} |_{X}$, $\cB_{Y} \coloneqq \cR^{n} |_{[a, b]}$, $\Xi^{X} \coloneqq \lambda^{n} |_{\cB_{X}}$ and $\Xi^{Y} \coloneqq \lambda^{n} |_{\cB_{Y}}$. Clearly, $f = g \circ h$, $\cB_{Y} \subseteq h^{\to}(\cB_{X})$ and $\Xi^{Y} \subseteq \Xi_{h}^{X}$, that is, we are under the hypothesis of~\autoref{da1}, by virtue of which it follows that $g \circ h \in \cI(\Xi^{X})$, or equivalently by~\autoref{t100}, $[f \in \cI(\lambda^{n} |_{[a, b]})]^{M}$. Moreover, by~(\ref{e66}) we have that
    $$
    \int_{[a, b]} g d \lambda^{n} |_{[a, b]} = \int_{[a, b]^{M}} f d \lambda^{n} |_{[a, b]^{M}}, 
    $$
    and therefore, by~\autoref{t100} again, it follows that:
    $$
    \left( \int_{[a, b]} f d \lambda^{n} |_{[a, b]} \right)^{M} = \left(  \int_{[a, b]} g d \lambda^{n} |_{[a, b]} \right)^{N}. 
    $$

    \ref{s20.1} $\Rightarrow$ \ref{s20.2}:  Intuitively, the situation in this case is illustrated in~\autoref{fig:enter-label}. When we consider that $f$ is Riemann integrable on $[a, b]$ in $M$, it follows that $f$ is continuous except on a Lebesgue measure zero set. However, when interpreting $f$ in the model $N$, the function $f$ is not defined on the entire rectangle $[a, b]$, as new real numbers may appear. Therefore, our proof consists of extending $f$ to these new real numbers in such a way that its integrability is preserved. To achieve this, we will approximate $f$ using step functions as follows. 

    Working in $M$, assume that $f$ is $\lambda^{n}|_{[a, b]}$-integrable. For any $m \in \bbN$, we can find step functions $\sigma_{m}, \tau_{m}$ on $[a, b]$ such that $\sigma_{m} \leq f \leq \tau_{m}$, and 
    \begin{equation}\label{e24}
        \text{$
    \int_{[a, b]} (\tau_{m} - \sigma_{m}) d \lambda^{n} |_{[a, b]} < \frac{1}{m + 1}.$}
    \end{equation}

    Since each $\sigma_{m}$ and $\tau_{m}$ are step functions, for any $m \in \bbN$ there are a partition $P_{m}$ of $[a, b]$ and sequences of real numbers $\bar{\alpha}_{m} = \langle \alpha_{m, R} \colon R \in \calS_{P_{m}} \rangle$, $\bar{\beta}_{m} = \langle \beta_{i, R} \colon R \in \calS_{P_{m}} \rangle$,  such that, without loss of generality, for any $R \in \calS_{P_{m}}$ and $x \in R^{\circ}$, $\sigma_{m}(x) = \alpha_{m, R}$ and $\tau_{m}(x) = \beta_{m, R}$, where $R^{\circ}$ denotes the interior of $R$.  

    Now, we work in $N$. For any $m \in \bbN$, let $\sigma_{m}^{\ast}, \tau_{m}^{\ast} \colon [a, b] \to \bbR$ be such that for any $x \in [a, b] \setminus \bigcup_{R \in \calS_{P_{m}}}R^{\circ}$, $\sigma_{m}^{\ast}(x) \coloneqq \alpha_{m, R_{0}}$ and $\tau_{m}^{\ast}(x) \coloneqq \beta_{m, R_{0}}$, where $R_{0}$ is the unique sub-rectangle determined by $P_{m}$ such that $x \in R_{0}^{\circ}$. If $x \in R^{\circ}$ for some $R \in \calS_{P_{m}}$, we can define  $\sigma_{m}^{\ast}(x) $ and $ \tau_{m}^{\ast}(x)$ arbitrarily. 
    Notice that this makes sense because we can use the end-points of the partitions to extend $P_{m}$ to a partition of $[a, b]$ in $N$. Clearly, $\sigma_{m}^{\ast}$, $\tau_{m}^{\ast}$ are step functions such that $\sigma_{m}^{\ast} \leq \tau_{m}^{\ast}$,  $\sigma_{m}^{\ast} \rest [a, b]^{M} \leq f \leq \tau_{m}^{\ast} \rest [a, b]^{M}$, and since the integral of step functions is merely a finite sum, for any $m \in \bbN$, 
    \begin{equation}\label{e34}
        \text{ $\int_{[a, b]} \tau^{\ast}_{m} d \lambda^{n} |_{[a, b]}  = \left( \int_{[a, b]} \tau_{m} d \lambda^{n} |_{[a, b]} \right)^{M}$ and \ $\int_{[a, b]} \sigma^{\ast}_{m} d \lambda^{n} |_{[a, b]}  = \left( \int_{[a, b]} \sigma_{m} d \lambda^{n} |_{[a, b]} \right)^{M}\! \! \! \!.$} 
    \end{equation}

    Define $h \colon [a, b] \to \bbR$ by $h(x) \coloneqq \inf_{m \in \bbN} \tau_{m}^{\ast}(x)$ whenever $x \in [a, b]$.  Using this, we can introduce the desired extension of $f$: define $g \colon [a, b] \to \bbR$, as follows: 
    $$
    g(x) \coloneqq 
    \begin{cases}
    f(x) \text{ if } x \in [a, b]^{M}, \\
      \\
    h(x) \text{ if } x \in [a, b] \setminus [a, b]^{M}.  
    \end{cases}
    $$

    Let $\varp > 0$ and $K \in \bbN$ such that $\frac{1}{K + 1} < \varp$. It is not hard to check that $\sigma_{K}^{\ast} \leq g \leq \tau^{\ast}_{K}$ except of a Lebesgue measure zero set, and by~(\ref{e24}) and~(\ref{e34}), it follows that 
    $$
    \int_{[a, b]} (\tau_{K}^{\ast} - \sigma_{K}^{\ast}) d \lambda^{n} |_{[a, b]} < \frac{1}{K + 1} < \varp. 
    $$
    Thus, $g$ is a Riemann integrable function extending $f$. Notice that, in this case~(\ref{e-s20.1}) follows from the proof of~\ref{s20.2} $\Rightarrow$~\ref{s20.2}. 

    Finally, we deal with the uniqueness of $g$. Assume that, in $N$, $g^{\ast}$ is a Riemann-integrable function on $[a, b]$ extending $f$. Consider $E_{g}$ and $E_{g^{\ast}}$ as the set of discontinuities of $g$ and $g^{\ast}$, respectively. Set $E \coloneqq E_{g} \cup E_{g^{\ast}}$, whose Lebesgue measure is zero. Let $\varp > 0$ and $x \in [a, b] \setminus E$. Pick some sequence points with rational coordinates  $\langle x_{m} \colon m \in \bbN \rangle$ in $[a, b]$ converging to $x$. Since $g$  and $g'$ are continuous at $x$, it follows that:
    $$
    g^{\ast}(x) = g^{\ast} \left( \lim_{m \to \infty} x_{m} \right) = \lim_{m \to \infty} g^{\ast}(x_{m}) = \lim_{m \to \infty} g(x_{m}) = g \left( \lim_{m \to \infty} x_{m} \right) = g(x). 
    $$

    Thus, $g = g^{\ast}$ on $[a, b] \setminus E$. 
\end{PROOF}

By applying~\autoref{s20} to $N = V$, we get the following result: 

\begin{corollary}\label{s20.c}
    Let $M$ be a transitive  model of $\ZFC$, $n \in \bbN$, and  $a, b \in \bbR^{n} \cap M$ with $a \leq b$. In $M$, let $f$ be a real-valued function on $[a, b]$. Then, the following statements are equivalent: 

    \begin{enumerate}[label = \normalfont (\roman*)]
        \item\label{s20.1.c} $f$ is Riemann integrable in $M$.  

        \item\label{s20.2.c} There exists some Riemann integrable function $g \colon [a, b] \to \bbR$ extending $f$. 
    \end{enumerate}

    If either~\ref{s20.1.c} or~\ref{s20.2.c} holds, then: 
    \begin{equation}\label{e-s20.1.c}
        \text{$  \int_{[a, b]} g d \lambda^{n} |_{[a, b]} = 
    \left( \int_{[a, b]} f d \lambda^{n} |_{[a, b]} \right)^{M}   \! \! \!.
    $}
    \end{equation}

    Moreover, the function $g$  in~\ref{s20.2.c} is unique except in a Lebesgue measure zero set: if $g^{\ast}$ is another Riemann integrable function on $[a, b]$ extending $f$, then there exists some Lebesgue measure zero set $E \subseteq [a, b]$ such that, for any $x \in [a, b] \setminus E$, $g^{\ast}(x) = g(x)$. 
\end{corollary}

There is an alternative proof of~\ref{s20.2} $\Rightarrow$~\ref{s20.1} in~\autoref{s20} using once approximations with step functions. In~\autoref{da20} below, we outline a sketch of this proof.

\begin{remark}\label{da20}
    Assume the same hypothesis as in~\autoref{s20}. We say that a step function $\sigma$ on $[a, b]$ is \emph{rational} if its constant values and the endpoints of the partitions defining it —except possibly for the end-points $a_{i}$ and $b_{i}$ for $i < n$— are rational numbers. Observe that, for any step function $\sigma$ and any $\varepsilon > 0$, we can construct rational step functions $\sigma_{\varp, -}$ and $\sigma_{\varp, +}$ such that $\sigma_{\varp, -} \leq \sigma \leq \sigma_{\varp, +}$,  and  
    \begin{equation*}
        \text{ $\int_{[a, b]} \sigma d \lambda^{n} |_{[a, b]} - \varp \leq \int_{[a, b]} \sigma_{\varp, -} d\lambda^{n} |_{[a, b]} \text{  \ and  } \int_{[a, b]} \sigma d\lambda^{n} |_{[a, b]} + \varp \geq \int_{[a, b]} \sigma_{\varp, +} d\lambda^{n} |_{[a, b]}$.}
    \end{equation*}
    
    Working in $N$, assume that $g$ is a Riemann integrable function extending $f$ and let $\varp > 0$. We can find step functions $\sigma, \tau \colon [a, b] \to \bbR$ such that $\sigma \leq g \leq \tau$ and 
    $$
    \int_{[a, b]} (\tau - \sigma) d \lambda^{n} |_{[a, b]} < \frac{\varp}{2}. 
    $$
    
    Working in $M$ now, consider $\sigma_{\epsilon, -}$ and $\tau_{\epsilon, +}$ rational step functions as above, where $\epsilon \coloneqq \frac{\varp}{4}$. Since rational numbers are absolute for transitive models of $\ZFC$, we have that, $\sigma_{\epsilon, -}$ and $\tau_{\epsilon, +}$ are step functions such that $\sigma_{\epsilon, -} \leq f \leq \tau_{\epsilon, +}$ and 
    $$
    \int_{[a, b]} (\tau_{\epsilon, +} - \sigma_{\epsilon, -}) d \lambda^{n} |_{[a, b]} < \varp.
    $$
    Thus, $f$ is Riemann integrable in $M$. 

    We now deal with the value of the integrals. Without loss of generality, we can assume that $f$ and $g$ are non-negative functions. On the one hand, working in $N$, let $\rho$ be a step non-negative step function on $[a, b]$ such that $\rho \leq g$ and $\varp > 0$. Then, 
    \begin{equation*}
        \text{$\int_{[a, b]} g d \lambda^{n} |_{[a, b]} \geq \int_{[a, b]} \rho d \lambda^{n} |_{[a, b]}$} 
    \end{equation*}

    Now, in $M$ we have that $\rho_{\varp, -}$ is step function such that $\rho_{\varp, -} \leq f$, and therefore,
    $$
    \left( \int_{[a, b]} f d \lambda^{n} |_{[a, b]} \right)^{M} \geq \left( \int_{[a, b]} \rho_{\varp, -} d \lambda^{n} |_{[a, b]} \right)^{M} = \int_{[a, b]} \rho_{\varp, -} d \lambda^{n} |_{[a, b]}  \geq \int_{[a, b]} \rho d \lambda^{n} |_{[a, b]} - \varp.
    $$
    Finally, since $\varp$ is arbitrary, it follows that 
    $$
    \left( \int_{[a, b]} f d \lambda^{n} |_{[a, b]} \right)^{M} \geq \left( \int_{[a, b]} g d \lambda^{n} |_{[a, b]} \right)^{N}. \! \! \! 
    $$ 
    
    The proof for the converse inequality is similar by considering the extension of step functions as in the proof of \ref{s20.1} $\Rightarrow$~\ref{s20.2}  in~\autoref{s20}.  
\end{remark}

We complete this paper by stating some natural questions that arose from~\autoref{s20}. 

\begin{question}\label{q100}\ 

    \begin{enumerate}[label = \normalfont (\arabic*)]
        \item\label{q100.2} Is the Riemann-Stieltjes integral absolute for transitive models of $\ZFC$? 
        
        \item\label{q100.3} Is the Lebesgue integral absolute for transitive models of $\ZFC$?
    \end{enumerate}
\end{question}

\autoref{q100}~\ref{q100.3} is particularly interesting. Although one approach to defining the Lebesgue integral involves simple functions, the methods used to prove~\autoref{s20} do not apply in this case. For instance, if $M \subseteq N$ are transitive models of $\ZFC$ and $N$ contains a Cohen real over $M$, then $[a, b]^{M}$ has Lebesgue measure zero in $N$. This makes the information provided by $f$ completely irrelevant —in the context of the Lebesgue integral— for defining a potential function $g$ as in~\autoref{s20}, because there are functions that are Lebesgue integrable and discontinuous everywhere. This suggests that the uniqueness achieved previously for the Riemann integral  probably does not hold for the Lebesgue integral, since for this the characterization of integrability in terms of continuity, provided by the Lebesgue-Vitali theorem, is absolutely fundamental.

{\small
\bibliographystyle{appl}

\bibliographystyle{alpha}
}

\end{document}